\documentclass[11pt,reqno]{amsart}

\ifdefined\synctex
\fi

\usepackage[T1]{fontenc}
\usepackage{lmodern}
\usepackage{amsmath,amssymb,mathtools}
\usepackage[hmargin=1.08in,vmargin=1.00in]{geometry}
\usepackage[expansion=false]{microtype}
\usepackage[
  bookmarksnumbered,
  colorlinks=true,
  linkcolor=black,
  citecolor=black,
  urlcolor=blue
]{hyperref}
\hypersetup{
  pdftitle={The regular polygon maximizes the first positive Steklov
    eigenvalue at fixed area},
  pdfsubject={Spectral optimization of fixed-area polygons}
}

\numberwithin{equation}{section}
\newtheorem{theorem}{Theorem}[section]
\newtheorem{lemma}[theorem]{Lemma}
\newtheorem{numericallemma}[theorem]{Numerical Lemma}
\newcommand{\R}{\mathbb R}
\allowdisplaybreaks[1]
\ifdefined\pdfoutput
  \ifnum\pdfoutput>0
    \IfFileExists{embedfile.sty}{%
      \usepackage{embedfile}%
      \AtBeginDocument{\embedfile[mimetype=text/plain,desc={LaTeX source}]{\jobname.tex}}%
    }{}%
  \fi
\fi

\title[The fixed-area polygonal Steklov inequality]{The regular polygon
  maximizes the first positive Steklov eigenvalue among simple polygons at fixed area}
\author[S.~Luo]{Senping Luo}
\author[B.~Xie]{Bingbing Xie}

\address[S.~Luo and B.~Xie]{School of Mathematics and statistics,
  Jiangxi Normal University, Nanchang, 330022, China; and Jiangxi applied
  mathematical research center.}
\email[S.~Luo]{luosp1989@163.com }
\email[B.~Xie]{xiebb2074@163.com }

\begin{document}
\raggedbottom

\begin{abstract}
  For $N\ge5$ and every simple planar polygon $P$ with at most $N$ sides,
  $\sqrt{|P|}\,\sigma_1(P)\le\sqrt{|R_N|}\,\sigma_1(R_N)$, where $\sigma_1$ is the
  first positive Steklov eigenvalue and $R_N$ is a regular $N$-gon.
  Equality holds precisely for polygons similar to $R_N$.
  This answers Oudet's fixed-area open problem \cite{AIM} for all $N\ge5$.

\end{abstract}
\maketitle

\section{Introduction}
Let $P\subset\R^2$ be the interior of a polygonal Jordan curve, with area
$|P|=\int_Pdx$ and perimeter $L(P)=\int_{\partial P}ds$.
The first positive Steklov eigenvalue $\sigma_1(P)$ is the least $\sigma>0$ for which
\begin{equation}
  \label{pde}
  \Delta u=0\quad\hbox{in }P,\qquad
  \partial_\nu u=\sigma u\quad\hbox{on }\partial P
\end{equation}
has a nonzero weak solution $u\in H^1(P)$, that is,
$\int_P\nabla u\cdot\nabla v\,dx=\sigma\int_{\partial P}uv\,ds$ for $v\in H^1(P)$.
Here $\nu$ is the outward unit normal and $\partial_\nu u=\nabla u\cdot\nu$ on each
open side.
Let $R_N$ be a regular $N$-gon.
Oudet \cite[Problem~4.6]{AIM} asked whether $R_N$ maximizes $\sigma_1$ among $N$-gons
of prescribed area or perimeter.
Since $\sigma_1(tP)=t^{-1}\sigma_1(P)$ for $t>0$, the fixed-area question amounts to
maximizing $\sqrt{|P|}\,\sigma_1(P)$.

\begin{theorem}
  \label{main}
  Let $N\ge5$.
  Every simple planar polygon $P$ with at most $N$ sides satisfies
  \begin{equation}
    \label{mainineq}
    \sqrt{|P|}\,\sigma_1(P)\le\sqrt{|R_N|}\,\sigma_1(R_N).
  \end{equation}
  Equality holds if and only if $P$ is similar to $R_N$.
\end{theorem}

Equivalently, among simple polygons with at most $N$ sides and the same area, $R_N$
minimizes the sharp boundary trace constant

\begin{equation}
  \label{ray}
  \frac{1}{\sigma_1(P)}
  =\sup_{\substack{u\in H^1(P)\\\int_P|\nabla u|^2\,dx>0}}
  \frac{\displaystyle\int_{\partial P}\!
    \left|u(x)-\frac{1}{L(P)}\int_{\partial P}u(y)\,ds_y\right|^2ds_x}
    {\displaystyle\int_P|\nabla u|^2\,dx}.
\end{equation}

Weinstock \cite{W} and Brock \cite{Brock} proved the sharp disk bounds
\[
  L(P)\sigma_1(P)\le 2\pi,\qquad \sqrt{|P|}\,\sigma_1(P)\le\sqrt\pi .
\]
Girouard and Polterovich give both inequalities in Henrot's volume \cite[Chapter~5,
Section~5.4]{Henrot}.
Freitas and Laugesen \cite{FL} recovered Weinstock's bound from the second Robin
eigenvalue.
For a disk $D$ of area $|P|$, Brasco, De Philippis and Ruffini
\cite[Corollary~5.2]{BDR} proved
\[
  \sqrt\pi-\sqrt{|P|}\,\sigma_1(P)\ge
  c\inf_{y\in\mathbb R^2}
  \left(\frac{|P\mathbin{\triangle}(y+D)|}{|P|}\right)^2 ,
\]
where $c>0$ is absolute and $P\mathbin{\triangle}Q=(P\setminus Q)\cup(Q\setminus P)$.
The coordinate trials in the proof of Brasco and De Philippis \cite[Section~7.5.3,
Theorem~7.45]{Henrot} are $|P|^{-1/2}(x_i-L(P)^{-1}\int_{\partial P}x_i\,ds)$,
$i=1,2$.

For a polygon with at most $N$ sides, polygonal isoperimetry gives
$L(P)^2\ge4N\tan(\pi/N)|P|$, with equality precisely for polygons similar to $R_N$.
Indrei and Nurbekyan \cite{IN} quantified the deficit in terms of the side lengths,
vertex radii and deviation from convexity.
Weinstock's inequality therefore gives
\[
  \sqrt{|P|}\,\sigma_1(P)\le\frac{\pi}{\sqrt{N\tan(\pi/N)}},\qquad
  \sqrt{|R_N|}\,\sigma_1(R_N)
  =\frac{L(R_N)\sigma_1(R_N)}{2\sqrt{N\tan(\pi/N)}} .
\]
The second right-hand side is the first multiplied by $L(R_N)\sigma_1(R_N)/(2\pi)<1$.

Cheng, Gui, Hu, Li and Yao \cite[Theorems~1.1--1.3]{CGHLY} studied the first positive
Steklov eigenvalue of regular polygons at fixed perimeter.
Taking $L(R_N)=2\pi$ for every $N\ge3$, they proved
\[
  \begin{aligned}
    \sigma_1(R_N) &<\sigma_1(R_{N+1}) && (N\ge3), \\
    \sigma_1(R_N) &=1-\frac{2\zeta(3)}{N^3}-\frac{8\zeta(4)}{N^4}
      -\frac{26\zeta(5)}{N^5}+O(N^{-6}),
  \end{aligned}
\]
where $\zeta(s)=\sum_{j=1}^{\infty}j^{-s}$ for $s>1$.
They also proved that the real eigenspace for $\sigma_1(R_N)$ has dimension two.
With $L(R_N)=2\pi$, one has $|R_N|=\pi^2/(N\tan(\pi/N))$.
Rescaling to unit area and using $N\tan(\pi/N)\downarrow\pi$ gives
\[
  \begin{aligned}
    \sqrt{|R_N|}\,\sigma_1(R_N) &<\sqrt{|R_{N+1}|}\,\sigma_1(R_{N+1}) && (N\ge3), \\
    \sqrt{|R_N|}\,\sigma_1(R_N) &=\sqrt\pi\left(1-\frac{\pi^2}{6N^2}
      -\frac{2\zeta(3)}{N^3}+O(N^{-4})\right).
  \end{aligned}
\]
At area $A>0$ the maximum for $N\ge5$ is
$L(R_N)\sigma_1(R_N)/(2\sqrt{AN\tan(\pi/N)})$.
The infimum is zero.
Indeed, on $P_t=(0,t)\times(0,A/t)$ the trial $u(x,y)=x-t/2$ in \eqref{ray} gives
\[
  0<\sigma_1(P_t)\le\frac{A}{t^3/6+At/2}\longrightarrow0\qquad(t\to\infty).
\]

The corresponding fixed-perimeter functional is $L(P)\sigma_1(P)$, also invariant
under similarities.
Davies \cite[p.~1]{D} records both normalizations of Oudet's question.
Nigam \cite[Problem~8]{Nigam} posed the perimeter question for all $N$-gons.
Cheng, Gui, Hu, Li and Yao \cite[Conjecture~1.1]{CGHLY} conjectured that, among convex
polygons with at most $N$ sides and perimeter $2\pi$, the regular $N$-gon uniquely
maximizes $\sigma_1$ up to congruence.

For triangles, Dom\'inguez, Nigam and Shahriari \cite{DNS} gave numerical evidence for
the perimeter conjecture.
Davies \cite[Theorem~1.1]{D} gives the computer-assisted bound
\[
  L(R_3)\sigma_1(R_3)-L(P)\sigma_1(P)
  \ge\frac{\ell_{\max}-\ell_{\min}}{2\ell_{\max}},
\]
where $\ell_{\min}$ and $\ell_{\max}$ are the shortest and longest side lengths.
Girouard, Lagac\'e, Polterovich and Savo \cite[Theorem~1.6]{GLPS} proved the square
bound among rectangles.
Dong, Yang and Pan \cite[Theorem~1]{DYP} give
$L(P)\sigma_1(P)\le L(R_4)\sigma_1(R_4)=8t\tanh t$ for convex polygons with at most
four sides, where $\tan t\tanh t=1$ and $0<t<\pi/2$.
For the square spectrum by separation of variables, see \cite[Section~5.3.1]{Henrot}.
In both classes polygonal isoperimetry converts the perimeter bound to an area bound
through the identity
\[
  \frac{\sqrt{|P|}\,\sigma_1(P)}{\sqrt{|R_N|}\,\sigma_1(R_N)}
  =\frac{L(P)\sigma_1(P)}{L(R_N)\sigma_1(R_N)}
  \frac{\sqrt{|P|}/L(P)}{\sqrt{|R_N|}/L(R_N)}
  \le\frac{L(P)\sigma_1(P)}{L(R_N)\sigma_1(R_N)} .
\]

For affine images with $|\det T|=1$, Theorem~\ref{aff-main} also gives
\[
  \frac{L(TR_N)\sigma_1(TR_N)}{L(R_N)\sigma_1(R_N)}
  \le1-\frac{1-(\max_{|x|=1}|Tx|)^{-1}}{50}\qquad(N\ge5).
\]
For every equiangular pentagon $P$, Theorem~\ref{equiancrefine} gives
$L(P)\sigma_1(P)\le L(R_5)\sigma_1(R_5)$.

For the higher Steklov eigenvalues $0=\sigma_0(P)<\sigma_1(P)\le\sigma_2(P)\le\cdots$,
Hersch, Payne and Schiffer \cite{HPS} obtained $L(P)\sigma_k(P)\le2\pi k$ for $k\ge1$.
If $P$ is a convex $N$-gon with smallest interior angle $\alpha$, Dryden, Gordon,
Moreno, Rowlett and Villegas-Blas \cite[Theorem~3.10]{DGMRV} proved
$L(P)\sigma_k(P)\le800\pi^3(N-1)^2k^2\alpha/49$.
For curvilinear polygons $P,Q$ with smooth sides and the same cyclic sequence of side
lengths and corresponding angles in $(0,\pi)$, Levitin, Parnovski, Polterovich and
Sher \cite[Corollary~1.6]{LPPS} proved $\sigma_k(P)-\sigma_k(Q)=O(k^{-\varepsilon})$
as $k\to\infty$, for some $\varepsilon>0$.

At fixed area, P\'olya and Szeg\H{o} \cite{PS} proposed the Dirichlet inequality
$\lambda_1^D(P)\ge\lambda_1^D(R_N)$.
Dahne, G\'omez-Serrano and Pech-Alberich \cite[Theorem~1.1]{DGP} proved
$\lambda_1^D(R_{N+1})<\lambda_1^D(R_N)$ for the regular polygons of equal area,
$N\ge3$.
Cheng, Gui, Hu and Li \cite{CGHL} proved the polygonal inequality for convex polygons
with sufficiently many sides, Wang \cite{Wang} for simple polygons and $N\ge5$, and
Y.~Li \cite{YLi} for simple polygons and $N\ge3$.
For the torsional rigidity $T(P)=\int_P u\,dx$, where $-\Delta u=1$ in $P$ and $u=0$
on $\partial P$, Gui, Hu and Li \cite[Corollary~1.2]{GHL} proved
\[
  T(R_N)-T(P)\ge\frac{|P|^2}{8N\tan(\pi/N)}
  \left(1-\frac{L(R_N)^2}{L(P)^2}\right)^2
\]
when $P$ is an $N$-gon circumscribed about a circle and $|P|=|R_N|$.
Gui, Hu, Li and Zhang \cite{GHLZ} proved $T(R_{N+1})>T(R_N)$ at fixed area.
Gui, Hu, Li and Yang \cite{GHLY} obtained the corresponding tangential polygon
comparison under the Robin boundary condition $\partial_\nu u+\beta u=0$, $\beta>0$.
Wang \cite{Wang} and Bogosel, Bucur and Fragal\`a \cite{BBF} proved $T(P)\le T(R_N)$
for simple polygons of the same area.

\medskip
\noindent\textit{Main ideas of the proof.}
Laugesen and Siudeja \cite{LS} transplanted Neumann eigenfunctions from the
equilateral triangle.
In our work on the polygonal Neumann conjecture \cite{LX}, triangular transplantation
gives $|P|\mu_1(P)\le |R_N|\mu_1(R_N)$ for the first positive Neumann eigenvalue
$\mu_1$.
Laugesen and Siudeja recorded this conjecture in \cite[Conjecture~6.36]{Henrot}.
Here the constraint $\int_{\partial P}W\,ds=0$ dictates the choice of center.
By \eqref{eq:convex-reduction}, it suffices to consider a convex $N$-gon.
Translate and dilate so that $\int_{\partial P}x\,ds=0$ and $|P|=|R_N|$.
Let $v_j,p_j$ be the counterclockwise vertices of the centered $R_N$ and of $P$,
respectively, with $v_N=v_0$ and $p_N=p_0$.
Identify $z=x+iy$ with $(x,y)^{\mathsf T}$.
For $0\le j<N$, let
\begin{equation}
  \label{eq:sector-map-coordinates}
  \begin{gathered}
    S_j=\{sv_j+tv_{j+1}:s,t\ge0,\ s+t\le1\}, \\
    T_j=
    \begin{pmatrix}
      p_{j,1} & p_{j+1,1} \\
      p_{j,2} & p_{j+1,2}
    \end{pmatrix}
    \begin{pmatrix}
      v_{j,1} & v_{j+1,1} \\
      v_{j,2} & v_{j+1,2}
    \end{pmatrix}^{-1}.
  \end{gathered}
\end{equation}
The first eigenspace of $R_N$ has dimension two \cite[Theorem~1.2]{CGHLY}.
Section~\ref{sec:regular} chooses a pair $U=(U_1,U_2)$ satisfying $U(gx)=gU(x)$ for
every symmetry $g$ of $R_N$.
Define $W(T_jx)=T_jU(x)$ for $x\in S_j$.
Reflection across each common radius gives matching traces.
On the boundary, reflection across the perpendicular bisector of each side reduces
$\int_{\partial P}W\,ds=0$ to
\[
  \sum_{j=0}^{N-1}|p_{j+1}-p_j|(p_j+p_{j+1})=2\int_{\partial P}x\,ds=0.
\]
The shape of sector $j$ is determined up to similarity by the angle
$\arg(p_{j+1}/p_j)\in(0,\pi)$ and the ratio $|p_{j+1}|/|p_j|$.
Normalize $\int_{R_N}\sum_{k=1}^2|\nabla U_k|^2\,dx=2\pi$.
With the sector moments $b,B,C$ of \eqref{sect-moments}, Lemma~\ref{sect-scalar} gives
\[
  \begin{aligned}
    \frac{N}{2\pi}\bigg( & \int_{T_j(S_j)}\sum_{k=1}^2|\nabla W_k|^2\,dx
      -\sigma_1(R_N)|p_{j+1}-p_j|\int_0^1|W((1-t)p_j+tp_{j+1})|^2\,dt\bigg) \\
    & \le\frac{1-\det T_j}{2}
      +\frac{2(C-B)-3/2+2b}{\sin(2\pi/N)}
      \left(\arg(p_{j+1}/p_j)-\frac{2\pi}{N}\right).
  \end{aligned}
\]
The two terms on the right sum to zero over $0\le j<N$, since
\[
  \sum_{j=0}^{N-1}\det T_j=N,\qquad
  \sum_{j=0}^{N-1}\arg(p_{j+1}/p_j)=2\pi .
\]
Consequently
\[
  \sigma_1(P)\int_{\partial P}|W|^2\,ds
  \le\sum_{j=0}^{N-1}\int_{T_j(S_j)}\sum_{k=1}^2|\nabla W_k|^2\,dx
  \le\sigma_1(R_N)\int_{\partial P}|W|^2\,ds .
\]
Section~\ref{sec:comparison} proves the sector estimate by polynomial inequalities for
$5\le N\le314$ and a uniform bound for $N\ge315$.
Equality forces $p_{j+1}=e^{2\pi i/N}p_j$ for $0\le j<N$.
Appendix~\ref{app:numerical} gives the finite numerical estimates.

\section{The first eigenspace and affine images of the regular polygon}
\label{sec:regular}

Directional estimates for the first eigenspace of $R_N$ give a quantitative stability
bound for affine images at fixed area and at fixed perimeter.

\begin{theorem}[Affine stability]
  \label{aff-main}
  Let $N\ge5$ and $T\in\R^{2\times2}$ satisfy $|\det T|=1$.
  Then
  \begin{equation}
    \label{affresult}
    \frac{\sqrt{|TR_N|}\,\sigma_1(TR_N)}{\sqrt{|R_N|}\,\sigma_1(R_N)}
    \le\frac{L(TR_N)\sigma_1(TR_N)}{L(R_N)\sigma_1(R_N)}
    \le 1-\frac{1-(\max_{|x|=1}|Tx|)^{-1}}{50}.
  \end{equation}
  Among affine images of $R_N$, either functional is maximized precisely by its images
  under similarities.
\end{theorem}

\begin{lemma}
  \label{gap}
  Center $R_n$ at $0$ with a vertex on the positive real axis.
  In either rotation class $u(e^{2\pi i/n}z)=e^{\pm2\pi i/n}u(z)$, the second Steklov
  eigenvalue of $R_n$, multiplied by $L(R_n)/(2\pi)$, is at least
  \[
    \gamma_n=\frac{n(n-1)\sin(\pi/n)\cos(\pi/n)}
    {\pi\big(\sqrt{1+\tan^2(\pi/n)/4}+\tfrac12\tan(\pi/n)\big)^2}
    >\frac{37}{25}\qquad(n\ge5).
  \]
  For a unit vector $e=(e_1,e_2)$ let $e^\perp=(-e_2,e_1)$.
  For $U_n$ in \eqref{eq:regular-eigenfunction}, the real first eigenfunction
  \[
    u_e(x)=
    \frac{e_1\operatorname{Re}\big(U_n(2\pi x/L(R_n))\big)
    +e_2\operatorname{Im}\big(U_n(2\pi x/L(R_n))\big)}
    {\sqrt\pi}
  \]
  satisfies
  \begin{equation}
    \label{direction}
    \int_{R_n}|e_1\partial_{x_2}u_e-e_2\partial_{x_1}u_e|^2\,dx
    <\frac16\int_{R_n}\big(|\partial_{x_1}u_e|^2+|\partial_{x_2}u_e|^2\big)\,dx .
  \end{equation}
\end{lemma}

\begin{lemma}
  \label{uniformregular}
  For $n\ge3$ write $s_n=L(R_n)\sigma_1(R_n)/(2\pi)$, independent of the size of
  $R_n$.
  Then
  \[
    s_n>1-\frac{3.044}{n^3}\quad(n\ge20),\qquad
    s_n>1-\frac{2.435}{n^3}\quad(n\ge300).
  \]
\end{lemma}

For $n\ge5$ take $L(R_n)=2\pi$.
Then
\[
  R_n=\bigcap_{j=0}^{n-1}
    \left\{z\in\mathbb C:
    \operatorname{Re}\left(ze^{-(2j+1)\pi i/n}\right)
    <\frac{\pi}{n\tan(\pi/n)}\right\},
  \qquad \sigma_1(R_n)=s_n .
\]
The centered Schwarz--Christoffel map from $\mathbb D=\{z\in\mathbb C:|z|<1\}$ onto
$R_n$ is
\begin{equation}
  \label{eq:SC-map}
  \Phi_n(z)=\frac{\Gamma(1-1/n)^2}{\Gamma(1-2/n)}
    \int_0^z\frac{d\zeta}{(1-\zeta^n)^{2/n}},
  \qquad \Phi_n(0)=0,\quad \Phi_n'(0)>0 .
\end{equation}
Here $\Gamma(s)=\int_0^\infty t^{s-1}e^{-t}\,dt$ for $s>0$, and $(1-\zeta^n)^{-2/n}=1$
at $\zeta=0$ fixes the holomorphic branch \cite[Section~4.1]{DrTr}.
Euler's beta and reflection formulas give
\[
  \int_0^1(1-t^n)^{-2/n}\,dt
  =\frac{\Gamma(1+1/n)\Gamma(1-2/n)}{\Gamma(1-1/n)},\qquad
  \Phi_n(e^{2\pi ij/n})=\frac{\pi e^{2\pi ij/n}}{n\sin(\pi/n)}
  \quad(0\le j<n).
\]
Under $z=\Phi_n(e^{it})$, the boundary arclength is $ds=|\Phi_n'(e^{it})|\,dt$ almost
everywhere \cite[Proposition~2.4]{CGHLY}.
Explicitly,
\[
  |\Phi_n'(e^{it})|=\frac{\Gamma(1-1/n)^2}{\Gamma(1-2/n)}
    |2\sin(nt/2)|^{-2/n}.
\]
The Fourier coefficients vanish outside $n\mathbb Z$.
At the remaining frequencies Cheng, Gui, Hu, Li and Yao \cite[Lemma~3.1]{CGHLY} give
\[
  \rho_{n,|j|}:=\frac{1}{2\pi}\int_0^{2\pi}|\Phi_n'(e^{it})|e^{-injt}\,dt
  =\frac{\Gamma(|j|+1/n)\Gamma(1-1/n)}
    {\Gamma(1/n)\Gamma(|j|+1-1/n)}
  \qquad (j\in\mathbb Z).
\]
In particular, $\rho_{n,0}=1$ and $(nj+n-1)\rho_{n,j+1}=(nj+1)\rho_{n,j}$ for $j\ge0$.

Their first-eigenspace theorem \cite[Theorem~1.2]{CGHLY} identifies two conjugate
complex lines on which $u(e^{2\pi i/n}z)=e^{\pm2\pi i/n}u(z)$.
For a harmonic $u\in H^1(R_n)$ in the $+1$ rotation class let
\[
  h_j=\frac{\sqrt{|1+nj|}}{2\pi}
    \int_0^{2\pi}u(\Phi_n(e^{it}))e^{-i(1+nj)t}\,dt
    \qquad(j\in\mathbb Z).
\]
The harmonic expansion and the Dirichlet energy are
\[
  u(\Phi_n(re^{it}))=
  \sum_{j\in\mathbb Z}\frac{h_j}{\sqrt{|1+nj|}}\,
    r^{|1+nj|}e^{i(1+nj)t},
  \qquad \int_{R_n}\sum_{m=1}^{2}|\partial_{x_m}u|^2\,dx
    =2\pi\sum_{j\in\mathbb Z}|h_j|^2 .
\]
The boundary form of \cite[Proposition~3.3]{CGHLY} gives
\[
  \begin{aligned}
    \frac{1}{2\pi}\int_{\partial R_n}|u|^2\,ds
      &=\sum_{i,j\in\mathbb Z}A_{ij}h_j\overline{h_i}, &
      \frac{1}{s_n} &=\sup_{\sum_{j\in\mathbb Z}|h_j|^2=1}
      \sum_{i,j\in\mathbb Z}A_{ij}h_j\overline{h_i},
  \end{aligned}
\]
where
\[
  A_{ij}=\frac{\rho_{n,|i-j|}}{\sqrt{|1+ni|\,|1+nj|}}
  =\frac{\Gamma(|i-j|+1/n)\Gamma(1-1/n)}
  {\Gamma(1/n)\Gamma(|i-j|+1-1/n)\sqrt{|1+ni|\,|1+nj|}} .
\]
The matrix $A=(A_{ij})_{i,j\in\mathbb Z}$ defines a positive compact operator on
$\ell^2(\mathbb Z)$, since $\int_0^{2\pi}|\Phi'_n(e^{it})|^2\,dt<\infty$ for $n\ge5$
and $H^{1/2}(\partial\mathbb D)\hookrightarrow L^4(\partial\mathbb D)$ compactly.
Since $A_{00}=1$ and $A_{01}=1/((n-1)\sqrt{n+1})$,
\[
  \det
  \begin{pmatrix}
    A_{00}-1 & A_{01} \\
    A_{10} & A_{11}-1
  \end{pmatrix}
  =-\frac{1}{(n-1)^2(n+1)}<0,\qquad s_n^{-1}=\lambda_{\max}(A)>1 .
\]

An equivariant first eigenfunction is
\begin{equation}
  \label{eq:regular-eigenfunction}
  \begin{aligned}
    \xi_j &=\lim_{k\to\infty}
      \frac{(A^k)_{j0}}{\big(\sum_{i\in\mathbb Z}|(A^k)_{i0}|^2\big)^{1/2}}
      \qquad(j\in\mathbb Z), \\
    U_n(\Phi_n(re^{it})) &=\sum_{j\in\mathbb Z}\frac{\xi_j}{\sqrt{|1+nj|}}\,
      r^{|1+nj|}e^{i(1+nj)t}.
  \end{aligned}
\end{equation}
Indeed, $A_{ij}>0$ gives a nonnegative maximizing vector, and the eigenvalue equation
gives $\xi_i=s_n\sum_{j\in\mathbb Z}A_{ij}\xi_j>0$ for every $i\in\mathbb Z$.
Simplicity of the top eigenvalue \cite[Lemma~4.5]{CGHLY} and the spectral theorem
yield the displayed limit.
Thus $U_n(e^{2\pi i/n}z)=e^{2\pi i/n}U_n(z)$, $U_n(\bar z)=\overline{U_n(z)}$, and
\[
  \int_{R_n}\big(|\partial_{x_1}U_n|^2+|\partial_{x_2}U_n|^2\big)\,dx=2\pi,\qquad
  \int_{\partial R_n}|U_n|^2\,ds=2\pi/s_n .
\]

\begin{proof}[Proof of Lemma~\ref{uniformregular}]
  Let $K=(A_{ij})_{i,j\in\mathbb Z\setminus\{0\}}$.
  Schur's test with weights $|1+ni|^{-1/2}$ gives
  \[
    \lambda_{\max}(K)\le
    \sup_{i\in\mathbb Z\setminus\{0\}}\sum_{j\in\mathbb Z\setminus\{0\}}
      \frac{\rho_{n,|i-j|}}{|1+nj|}.
  \]
  The decreasing rearrangements of the two factors are bounded termwise by
  $(1,\rho_{n,1},\rho_{n,1},\ldots)$ and $(1/(n-1),1/n,1/n,1/(2n),1/(2n),\ldots)$.
  Hence
  \[
    \lambda_{\max}(K)\le\frac{1}{n-1}+\frac2n\sum_{j=1}^{\infty}\frac{\rho_{n,j}}{j},
    \qquad
    \sum_{i\in\mathbb Z\setminus\{0\}} A_{i0}^2
    =\frac2n\sum_{j=1}^{\infty}\frac{\rho_{n,j}^2}{j(1-1/(nj)^2)} .
  \]
  For $j>M\ge1$,
  \[
    \frac{\rho_{n,j}}{\rho_{n,M}}
    \le\exp\!\left(-(1-2/n)\sum_{i=M}^{j-1}\frac1{i+1}\right)
    \le\left(\frac{M+1}{j+1}\right)^{1-2/n}
    \le(1+1/M)(M/j)^{1-2/n}.
  \]

  The integral test yields
  \[
    \begin{aligned}
      \lambda_{\max}(K) & \le k(n,M):=\frac{1}{n-1}
        +\frac2n\left(\sum_{j=1}^M\frac{\rho_{n,j}}{j}
        +\frac{(1+1/M)\rho_{n,M}}{1-2/n}\right), \\
      n^3\sum_{i\in\mathbb Z\setminus\{0\}} A_{i0}^2 & \le B(n,M):=n^2
        \left(\sum_{j=1}^M\frac{2\rho_{n,j}^2}{j(1-1/(nj)^2)}
        +\frac{(1+1/M)^2\rho_{n,M}^2}
        {(1-2/n)(1-1/(n^2(M+1)^2))}\right).
    \end{aligned}
  \]
  If $k(n,M)<1$, the Schur complement at $s_n^{-1}>1$ gives
  \[
    1-s_n<s_n^{-1}-1
    =\sum_{\ell=0}^{\infty}s_n^{\ell+1}
      \sum_{\substack{i,j\in\mathbb Z\\ij\ne0}} A_{i0}(K^\ell)_{ij} A_{j0}
    \le\frac{B(n,M)}{n^3(1-k(n,M))} .
  \]
  Both bounds decrease with $n$, since
  $n\rho_{n,j}=\dfrac{n}{n-1}\displaystyle\prod_{i=1}^{j-1}\dfrac{ni+1}{n(i+1)-1}$ for
  $j\ge1$.
  Direct substitution gives
  \[
    \frac{B(20,4)}{1-k(20,4)}<3.044,\qquad
    \frac{B(300,16)}{1-k(300,16)}<2.435 ,
  \]
  with both denominators positive.
\end{proof}

For $5\le n\le19$ use the Schur complement at $z>\lambda_{\max}(K)$.
With $\mathrm{Id}$ the identity on $\ell^2(\mathbb Z\setminus\{0\})$,
\begin{equation}
  \label{schur}
  \sum_{\substack{i,j\in\mathbb Z\\ij\ne0}}A_{i0}[(z\,\mathrm{Id}
    -K)^{-1}]_{ij}A_{j0}<z-1
  \quad\Longrightarrow\quad s_n^{-1}<z .
\end{equation}
By \eqref{eq:schur-recursion}, $V_j\ge10^{12}\rho_{n,j}$ and $C_i\ge V_{|i|}/|1+ni|$.
The tail estimate in the proof of Lemma~\ref{uniformregular} gives
\[
  \rho_{n,j}\le\frac{257 V_{256}}{256\cdot10^{12}}
    \left(\frac{256}{j}\right)^{1-2/n}\qquad(j>256).
\]
Since $|1\pm nj|\ge(n-1)j$ for $j\ge1$,
\[
  \begin{aligned}
    \lambda_{\max}(K) & \le k:=\frac{1}{n-1}
      +\frac{2}{n\cdot10^{12}}\left(\sum_{j=1}^{256}\frac{V_j}{j}
      +\frac{257 V_{256}}{256(1-2/n)}\right), \\
    \sum_{|i|>256}A_{i0}^2 & \le T:=
      \frac{1}{(1-2/n)(n-1)}\left(\frac{257 V_{256}}{256\cdot10^{12}}\right)^2 .
  \end{aligned}
\]
Truncation at $|i|=256$ and $2\sqrt{BT}\le B/64+64T$ give
\[
  \begin{aligned}
    \sum_{i\in\mathbb Z\setminus\{0\}} A_{i0}^2
      & \le B:=T+10^{-24}\sum_{0<|i|\le256}C_i V_{|i|}, \\
    \sum_{\substack{i,j\in\mathbb Z\\ij\ne0}} A_{i0}A_{ij}A_{j0}
      & \le H:=10^{-36}\sum_{\substack{0<|i|\le256\\0<|j|\le256}} C_iC_jV_{|i-j|}
      +k(B/64+65T).
  \end{aligned}
\]
For $z>k$, $0\le K\le k\,\mathrm{Id}$ gives
\[
  (z\,\mathrm{Id}-K)^{-1}\le z^{-1}\mathrm{Id}
    +\frac{K}{z(z-k)}.
\]
Set $y=1/z$.
By \eqref{schur}, $s_n>s_n^{\rm lo}$ for
\begin{equation}
  \label{eq:certified-regular-bound}
  s_n^{\rm lo}=\max\big\{ y\in10^{-9}\mathbb Z\cap[0.8,1]:
    ky<1,\ (1-y-By^2)(1-ky)>Hy^3\big\}.
\end{equation}

Let $\theta=\pi/n$ and $|e|=1$.
Rotation gives $\int_{\partial R_n}(e\cdot x)\,ds=0$ and
$\int_{\partial R_n}(e\cdot x)^2\,ds=\frac12\int_{\partial R_n}|x|^2\,ds$.
At apothem one, a side has mean squared radius
$(2\tan\theta)^{-1}\int_{-\tan\theta}^{\tan\theta}(1+t^2)\,dt=1+\tan^2\theta/3$.
Thus
\begin{equation}
  \label{linearray}
  q_n:=\frac{L(R_n)|R_n|}{2\pi\int_{\partial R_n}(e\cdot x)^2\,ds}
    =\frac{n\tan(\pi/n)}{\pi(1+\tan^2(\pi/n)/3)}\le1,
\end{equation}
since $\arctan x\ge3x/(3+x^2)$ for $x\ge0$.

\begin{proof}[Proof of Lemma~\ref{gap}]
  On a sector of the circumradius-one polygon, with polar angle $|\varphi|\le\theta$
  measured from its side normal, set
  \[
    F(r,\varphi)=(r\cos\theta/\cos\varphi,\varphi),\qquad
    DF=\frac{\cos\theta}{\cos\varphi}
    \begin{pmatrix}
      1 & \tan\varphi \\
      0 & 1
    \end{pmatrix}
  \]
  in orthonormal polar frames.
  Extending by rotation gives $F:\mathbb D\to R_n$.
  Its singular-value ratio and boundary arclength factor satisfy
  \[
    \frac{\max_{|v|=1}|DF\,v|}{\min_{|v|=1}|DF\,v|}\le
    (\sqrt{1+\tan^2\theta/4}+\tan\theta/2)^2,\qquad
    \frac{ds_{F(e^{i\varphi})}}{d\varphi}
    =\frac{\cos\theta}{\cos^2\varphi}\le\frac1{\cos\theta}.
  \]
  The disk eigenfunctions $r^{|k|}e^{ik\varphi}$ have eigenvalue $|k|$
  \cite[Example~1.3.1]{GP}.
  In either rotation class the second eigenvalue is
  $\min_{j\in\mathbb Z\setminus\{0\}}|\pm1+jn|=n-1$.
  Applying min--max under $F$ and multiplying by $L(R_n)/(2\pi)=\sin\theta/\theta$
  gives the stated bound.
  Each factor of $\gamma_n$ is nondecreasing with $n$, so $\gamma_n\ge\gamma_5>1.486$.

  Normalize $L(R_n)=2\pi$.
  For every symmetry $g$ of $R_n$, $DU_n(gx)=gDU_n(x)g^{\mathsf T}$.
  Rotation and reflection give
  \[
    \int_{R_n}DU_n\,dx=\frac{I_2}{2}\int_{R_n}\operatorname{tr}(DU_n)\,dx,
    \qquad \int_{R_n}DU_n(DU_n)^{\mathsf T}\,dx
    =\frac{I_2}{2}\int_{R_n}|DU_n|^2\,dx=\pi I_2,
  \]
  where $I_2$ is the $2\times2$ identity.
  Thus
  \[
    \int_{R_n}|\nabla u_e|^2\,dx=1,\qquad
    \int_{R_n}\nabla u_e\cdot\nabla u_{e^\perp}\,dx
    =\frac{1}{\sqrt{|R_n|}}\int_{R_n}e\cdot\nabla u_{e^\perp}\,dx=0 .
  \]

  The Dirichlet projection of $x\mapsto(e\cdot x)/\sqrt{|R_n|}$ onto the real first
  eigenspace is therefore
  \[
    x\longmapsto\frac{u_e(x)}{\sqrt{|R_n|}}
      \int_{R_n}\partial_e u_e(y)\,dy,\qquad
    \partial_e u_e=e\cdot\nabla u_e .
  \]
  This linear trial function has Dirichlet energy one and belongs to the sum of the
  two rotation classes.
  Its Steklov expansion gives
  \[
    \frac1{q_n}=\frac{1}{|R_n|}\int_{\partial R_n}(e\cdot x)^2\,ds
    \le\frac1{\gamma_n}
      +\left(\frac1{s_n}-\frac1{\gamma_n}\right)
      \frac1{|R_n|}\left|\int_{R_n}\partial_e u_e\,dx\right|^2 .
  \]
  Here $q_n\le1$ by \eqref{linearray}, and $s_n>47/50$ by Numerical Lemma~\ref{finite}
  and Lemma~\ref{uniformregular}, since
  $\frac{n\tan(\pi/n)}{(n-1)\tan(\pi/(n-1))}\ge\frac{\pi}{6\tan(\pi/6)}>(47/50)^2$ for
  $7\le n\le19$.
  By Cauchy--Schwarz,
  \[
    \begin{aligned}
      \int_{R_n}|\partial_{e^\perp}u_e|^2\,dx
        &=1-\int_{R_n}|\partial_e u_e|^2\,dx
        \le 1-\frac1{|R_n|}\left|\int_{R_n}\partial_e u_e\,dx\right|^2 \\
      & \le\frac{1-s_n/q_n}{1-s_n/\gamma_n}
        \le\frac{1-47/50}{1-(47/50)/(37/25)}<\frac16 .
    \end{aligned}\qedhere
  \]
\end{proof}

\begin{proof}[Proof of Theorem~\ref{aff-main}]
  Normalize $T=A/\sqrt{|\det A|}$ and $L(R_N)=2\pi$.
  For $t=\max_{|x|=1}|Tx|$, choose $|e|=1$ with
  \[
    T^{\mathsf T}Te=t^2e,\qquad T^{\mathsf T}Te^\perp=t^{-2}e^\perp .
  \]
  On side $j$ let $\tau=\tau_j$ be its unit tangent.
  Since $dy=dx$ and $ds_{Tx}=|T\tau|\,ds_x$,
  \begin{align*}
    \int_{TR_N}|\nabla(v\circ T^{-1})|^2\,dy
      &=t^{-2}\int_{R_N}|\partial_e v|^2\,dx+t^2\int_{R_N}|\partial_{e^\perp}v|^2\,dx,
        \\
    \min_{c\in\R}\int_{\partial R_N}|v-c|^2|T\tau|\,ds
      & \ge t^{-1}\min_{c\in\R}\int_{\partial R_N}|v-c|^2\,ds .
  \end{align*}
  Taking $v=e\cdot x$ and
  $v=(e_1\operatorname{Re}(U_N(x))+e_2\operatorname{Im}(U_N(x)))/\sqrt\pi$ in
  \eqref{ray}, we obtain from \eqref{linearray} and \eqref{direction}
  \[
    \begin{aligned}
      \frac{\sigma_1(TR_N)}{\sigma_1(R_N)}
        & \le t^{-1}\min\{q_N/s_N,1+(t^4-1)/6\}, \\
      \frac{L(TR_N)^2}{L(R_N)^2}
        &=\left(\frac1N\sum_{j=0}^{N-1}|T\tau_j|\right)^2
        \le\frac1N\sum_{j=0}^{N-1}|T\tau_j|^2=\frac{t^2+t^{-2}}2 .
    \end{aligned}
  \]
  Hence $q_N\le1$ and $s_N>47/50$ give
  \[
    \frac{L(TR_N)\sigma_1(TR_N)}{L(R_N)\sigma_1(R_N)}
    \le \sqrt{\frac{1+t^{-4}}2}
      \min\left\{\frac{50}{47},1+\frac{t^4-1}{6}\right\}
    \le
    \begin{cases}
      1-(1-t^{-4})/32, & 1\le t^4\le3/2, \\[2pt]
      (50/47)\sqrt{5/6}<49/50, & t^4\ge3/2 .
    \end{cases}
  \]
  For $1\le t^4\le3/2$ the first bound follows from
  \[
    1-\frac{1+t^{-4}}2\left(1+\frac{t^4-1}{6}\right)^2
    =\frac{1-t^{-4}}{2}\left(1-\frac{t^4+1}{3}-\frac{t^8-1}{36}\right)
    \ge\frac{1-t^{-4}}{16}.
  \]
  Finally $|TR_N|=|R_N|$, $L(TR_N)\ge L(R_N)$ by polygonal isoperimetry, and $t=1$
  exactly when $T^{\mathsf T}T=I_2$.
\end{proof}

\section{Convex reduction and transplantation}
\label{sec:transplantation}

It suffices to consider convex $N$-gons, for which transplantation reduces
\eqref{mainineq} to a sector inequality.
For pentagons we also prove the perimeter bound in the equiangular class and local
quadratic stability at fixed area.

\begin{lemma}[Convex reduction]
  \label{hullmain}
  Let $N\ge5$.
  For a simple polygon $P$ with $k\le N$ vertices, let $H$ be its closed convex hull
  and $m$ the number of its extreme points.
  Then
  \begin{equation}
    \label{hullresult}
    \sup_{P:\,m<N,\ k\le N}\sqrt{|P|}\,\sigma_1(P)<\sqrt{|R_N|}\,\sigma_1(R_N).
  \end{equation}
  In particular,
  \begin{equation}
    \label{eq:convex-reduction}
    \sqrt{|P|}\,\sigma_1(P)\ge\sqrt{|R_N|}\,\sigma_1(R_N)\ \Longrightarrow\ N\le m
      \le k\le N
    \ \Longrightarrow\ \partial P=\partial H .
  \end{equation}
\end{lemma}

\begin{lemma}[Transplantation]
  \label{sect-scalar}
  Let $P$ be a convex $n$-gon, $n\ge5$, normalized by \eqref{eq:normalized-fan}.
  In the coordinates \eqref{sect-matrix}, suppose the moments \eqref{sect-moments}
  satisfy, on each sector,
  \begin{equation}
    \label{sect-sign}
    \frac{4\big(\tfrac32+B(x^2-1)+C(y^2-1)+K(xy-1)\big)^3}
      {27y\big(bx+(1-b)y\big)^2}
    \le \frac12+\left(2(C-B)-\frac32+2b\right)
      \frac{\phi-2\pi/n}{\sin(2\pi/n)},
  \end{equation}
  with equality only at $(r,d)=(1,0)$.
  Then $\sqrt{|P|}\,\sigma_1(P)\le\sqrt{|R_n|}\,\sigma_1(R_n)$, with equality only for
  regular $P$.
\end{lemma}

\begin{lemma}
  \label{moment}
  With $P,H,m$ as in Lemma~\ref{hullmain}, let $|P|=1$.
  Suppose that $P$ is the interior of $H$ or that $\partial P$ is obtained from
  $\partial H$ by replacing the side from $A$ to $B$ with the broken line from $A$
  through $D\in H$ to $B$.
  Then
  \begin{equation}
    \label{momineq}
    \frac{2}{\sigma_1(P)}\ge \frac4{L(H)}+\frac{L(H)^3}{12m^2}
      +\left(\frac8{L(H)}-\frac{L(H)}2\right)(|H|-1)
      +\frac{3(|H|-1)^2}{L(H)}.
  \end{equation}
\end{lemma}
\begin{proof}
  The coordinate trials of Bucur, Ferone, Nitsch and Trombetti \cite[proof of
  Theorem~3.1]{BFNT} give, by \eqref{ray} and $|P|=1$,
  \[
    \sigma_1(P)\min_{c\in\R^2}\int_{\partial P}|x-c|^2\,ds
    =\sigma_1(P)\int_{\partial P}\left|x
      -\frac{1}{L(P)}\int_{\partial P}y\,ds_y\right|^2ds_x
    \le\sum_{k=1}^2\int_P|\nabla x_k|^2\,dx=2.
  \]
  Write $c_H=L(H)^{-1}\int_{\partial H}x\,ds$, and let $\ell_i,m_i,\nu_i$ denote the
  length, midpoint and outward unit normal of side $i$, $1\le i\le m$.
  Cauchy--Schwarz and $\sum_{i=1}^m\ell_i^3\ge L(H)^3/m^2$ give
  \[
    \begin{aligned}
      \int_{\partial H}|x-c_H|^2\,ds &=\sum_{i=1}^m\left(\ell_i|m_i-c_H|^2
        +\frac{\ell_i^3}{12}\right) \\
      & \ge\frac1{L(H)}\left(\sum_{i=1}^m \ell_i(m_i-c_H)\cdot\nu_i\right)^2
        +\frac{L(H)^3}{12m^2}
        =\frac{4|H|^2}{L(H)}+\frac{L(H)^3}{12m^2},
    \end{aligned}
  \]
  by $\int_{\partial H}(x-c_H)\cdot\nu\,ds=\int_H 2\,dx=2|H|$.
  For a replaced side, take $A=(0,0)$, $B=(B_1,0)$, $B_1>0$, and $D=(D_1,D_2)$ with
  $H\subset\{x_2\ge0\}$.
  On the broken line $\gamma$ from $A$ through $D$ to $B$, $ds\ge dx_1$.
  Its enclosed triangle has area $|H|-1$.
  Green's formula gives
  \[
    \int_{\gamma}|x-c|^2\,ds-\int_0^{B_1}|(t,0)-c|^2\,dt
    \ge 2\int_{\triangle ADB}(x_2-c_2)\,dx
    =2(|H|-1)(D_2/3-c_2)\ge-2(|H|-1)c_2 .
  \]
  At arclength $s$ from $A$, $x_2\le\min(s,L(H)-s)$, so
  $(c_H)_2\le L(H)^{-1}\int_0^{L(H)}\min(s,L(H)-s)\,ds=L(H)/4$.
  Hence
  \[
    \begin{aligned}
      \min_{c\in\R^2}\int_{\partial P}|x-c|^2\,ds
        & \ge\frac{4|H|^2}{L(H)}+\frac{L(H)^3}{12m^2}
        +\min_{c\in\R^2}\big(L(H)|c-c_H|^2-2(|H|-1)c_2\big) \\
      &=\frac{4|H|^2}{L(H)}+\frac{L(H)^3}{12m^2}-2(|H|-1)(c_H)_2
        -\frac{(|H|-1)^2}{L(H)} .
    \end{aligned}
  \]
  For $P$ equal to the interior of $H$, $|H|=1$.
\end{proof}

\begin{proof}[Proof of Lemma~\ref{hullmain}]
  For the vertices $a_0,\ldots,a_{k-1}$ of $P$,
  \[
    H=\left\{\sum_{j=0}^{k-1}\lambda_j a_j:\lambda_j\ge0\ (0\le j<k),\
      \sum_{j=0}^{k-1}\lambda_j=1\right\},\qquad 3\le m\le k\le N .
  \]
  Weinstock's inequality, $L(H)\le L(P)$ and polygonal isoperimetry give
  \begin{equation}
    \label{hullbound}
    \sqrt{|P|}\,\sigma_1(P)\le 2\pi\frac{\sqrt{|P|}}{L(P)}
    \le 2\pi\frac{\sqrt{|H|}}{L(H)}
    \le \frac{\pi}{\sqrt{m\tan(\pi/m)}} .
  \end{equation}
  An edge of $P$ joining two extreme points $A,B$ is an edge of $H$.
  Indeed, $H\cap(A+\R(B-A))$ is the segment from $A$ to $B$.
  On the complementary open arc of $\partial P$,
  $\operatorname{Im}(\overline{(B-A)}(x-A))$ is continuous and nonzero, hence has
  constant sign.

  For $m<N$, \eqref{hullbound} and $m\tan(\pi/m)\ge (N-1)\tan(\pi/(N-1))$ give
  \[
    \frac{\sqrt{|P|}\,\sigma_1(P)}{\sqrt{|R_N|}\,\sigma_1(R_N)}
    \le\frac{1}{s_N}\sqrt{\frac{N\tan(\pi/N)}{(N-1)\tan(\pi/(N-1))}}.
  \]
  The bound is $<1$ for $7\le N\le19$ by Numerical Lemma~\ref{finite}.
  For $N\ge20$, Euler's product gives
  \[
    \frac{\sin(2x)}{2x}
    =\prod_{j=1}^{\infty}\left(1-\frac{4x^2}{\pi^2j^2}\right)\le e^{-2x^2/3},
    \qquad \frac{d}{dx}\log\frac{\tan x}{x}
    =\frac{2}{\sin(2x)}-\frac1x\ge\frac{2x}{3}
    \quad(0<x<\pi/2).
  \]
  Integration on $[\pi/N,\pi/(N-1)]$ and Lemma~\ref{uniformregular} give
  \[
    \sqrt{\frac{N\tan(\pi/N)}{(N-1)\tan(\pi/(N-1))}}\le e^{-\pi^2/(3N^3)}<1
      -3.1/N^3<s_N .
  \]

  For $N=5,6$, normalize $|P|=1$.
  Numerical Lemma~\ref{finite} gives
  \[
    \frac{77}{50}<\sqrt{|R_5|}\,\sigma_1(R_5),\qquad
    \frac{18}{11}<\sqrt{|R_6|}\,\sigma_1(R_6).
  \]
  Suppose $m<N$.
  By \eqref{hullbound},
  \[
    m\le3\ \Longrightarrow\ \sigma_1(P)
      \le\frac{\pi}{\sqrt{3\tan(\pi/3)}}<\frac{77}{50},
    \qquad m\le4\ \Longrightarrow\ \sigma_1(P)\le\frac{\pi}{2}<\frac{18}{11}.
  \]
  It remains to consider $m=N-1$.
  Since $k-m\le1$, Lemma~\ref{moment} applies.

  If $N=5$ and $\sigma_1(P)\ge77/50$, then
  \[
    \begin{gathered}
      L(H)\le L(P)\le100\pi/77<4.1, \\
      0\le |H|-1\le L(H)^2/16-1<1/16,\qquad 8/L(H)-L(H)/2>-1/10 .
    \end{gathered}
  \]
  Consequently
  \[
    \frac{2}{\sigma_1(P)}\ge\min_{x>0}\left(\frac4x+\frac{x^3}{192}\right)-\frac1{160}
    =\frac43-\frac1{160}>\frac{100}{77},
  \]
  a contradiction.
  If $N=6$ and $\sigma_1(P)\ge18/11$, then $L(H)\le L(P)\le11\pi/9<31/8<4$ and
  \[
    \frac{2}{\sigma_1(P)}\ge\frac4{L(H)}+\frac{L(H)^3}{300}
    \ge\frac{32}{31}+\frac{(31/8)^3}{300}>\frac{11}{9},
  \]
  since $-4/x^2+x^2/100<0$ for $0<x\le31/8$.
  By \eqref{hullresult}, $\sqrt{|P|}\,\sigma_1(P)\ge\sqrt{|R_N|}\,\sigma_1(R_N)$
  forces $k=m=N$, hence every edge lies on $\partial H$.
\end{proof}

For a convex $n$-gon $P$ with vertices $a_j$ (indices modulo $n$) and $|\nu|=1$,
\[
  \max_{x\in\overline P}\nu\cdot x
  -\frac{\nu}{L(P)}\cdot\int_{\partial P}x\,ds
  =\frac1{L(P)}\sum_{j=0}^{n-1}|a_{j+1}-a_j|
  \left(\max_{0\le k<n}\nu\cdot a_k-\nu\cdot\frac{a_j+a_{j+1}}2\right)>0 ,
\]
since at least one side lies off the supporting line.
Thus $L(P)^{-1}\int_{\partial P}x\,ds\in P$.
Order the vertices counterclockwise.
For $\theta=\pi/n$ and $R_n$ of circumradius one,
\[
  v_j=e^{2ij\theta},\qquad R_n=\bigcap_{j=0}^{n-1}
  \{z\in\mathbb C:\operatorname{Re}(e^{-(2j+1)i\theta}z)<\cos\theta\}.
\]
Translating the boundary centroid to $0$ and scaling by $\sqrt{|R_n|/|P|}$ give the
vertices
\begin{equation}
  \label{eq:normalized-fan}
  p_j=\sqrt{\frac{n\sin(2\theta)}{2|P|}}
  \left(a_j-\frac{1}{2L(P)}\sum_{i=0}^{n-1}
  |a_{i+1}-a_i|(a_i+a_{i+1})\right),\qquad 0\le j<n .
\end{equation}
Denote the normalized polygon again by $P$.
Then $|P|=|R_n|$, $\int_{\partial P}x\,ds=0$ and $0\in P$.
Hence $0<\phi_j:=\arg(p_{j+1}/p_j)<\pi$ for $0\le j<n$ and
$\sum_{j=0}^{n-1}\phi_j=2\pi$.
By \eqref{eq:sector-map-coordinates},
\[
  \begin{gathered}
    \overline P=\bigcup_{j=0}^{n-1}T_j(S_j),\qquad
      T_j(S_j)=\{sp_j+tp_{j+1}:s,t\ge0,\ s+t\le1\}, \\
    \det T_j=\frac{\operatorname{Im}(\overline{p_j}p_{j+1})}{\sin(2\theta)}>0,\qquad
      \sum_{j=0}^{n-1}\det T_j=\frac{2|P|}{\sin(2\theta)}=n .
  \end{gathered}
\]
On $\Gamma_j=\{(1-t)v_j+tv_{j+1}:0\le t\le1\}$ the outward normal and unit tangent are
$\nu_j=e^{(2j+1)i\theta}$ and $\tau_j=i\nu_j$.
On this sector write
\begin{equation}
  \label{sect-matrix}
  \begin{gathered}
    x=\frac{|T_j\nu_j|^2}{\det T_j},\qquad y=\frac{|T_j\tau_j|^2}{\det T_j},\qquad
      \phi=\arccos\frac{p_j\cdot p_{j+1}}{|p_j|\,|p_{j+1}|}, \\
    t=\tan(\pi/n),\qquad r=\frac{\tan(\phi/2)}{t},\qquad
      d=\frac{(|p_{j+1}|-|p_j|)^2}{|p_{j+1}|\,|p_j|}.
  \end{gathered}
\end{equation}
From $T_jv_j=p_j$, $T_jv_{j+1}=p_{j+1}$ and $0<\phi<\pi$,
\[
  x=\frac{1+d(1+t^2r^2)/4}{r},\qquad
  y=r+\frac{d(1+t^2r^2)}{4t^2r},\qquad
  \frac{(T_j\nu_j\cdot T_j\tau_j)^2}{(\det T_j)^2}=xy-1 .
\]
Rescaling \eqref{eq:regular-eigenfunction} gives
\[
  U(\xi)=U_n\!\left(\frac{\pi\xi}{n\sin\theta}\right),\qquad
  U(e^{2i\theta}\xi)=e^{2i\theta}U(\xi),\qquad
  U(\bar\xi)=\overline{U(\xi)} .
\]
Write $u=U\cdot\nu_j$, $v=U\cdot\tau_j$ and $\partial_e=e\cdot\nabla$.
By rotation, the integrals are independent of $j$.
\begin{equation}
  \label{sect-moments}
  \begin{gathered}
    \int_{S_j}|DU|^2\,d\xi=\frac{2\pi}{n},\qquad
      \int_{\Gamma_j}|U|^2\,ds=\frac{2\sin\theta}{s_n}, \\
    b=\frac{s_n}{2\sin\theta}\int_{\Gamma_j}u^2\,ds,\quad
      B=\frac{n}{2\pi}\int_{S_j}(\partial_{\tau_j}u)^2\,d\xi,\quad
      C=\frac{n}{2\pi}\int_{S_j}(\partial_{\nu_j}v)^2\,d\xi, \\
    K=\frac{n}{2\pi}\int_{S_j}\big[(\partial_{\nu_j}u-\partial_{\tau_j}v)^2
      -2(\partial_{\tau_j}u)(\partial_{\nu_j}v)\big]\,d\xi .
  \end{gathered}
\end{equation}
Set $W(T_j\xi)=T_jU(\xi)$ for $\xi\in S_j$.
Reflection in $\R v_j$ gives, for $0<s<1$,
\[
  U(sv_j)\in\R v_j,\qquad T_{j-1}U(sv_j)=T_jU(sv_j).
\]
Thus $W\in H^1(P;\R^2)$.
Reflection in $\R\nu_j$ gives $u(\nu_j^2\bar\xi)=u(\xi)$ and
$v(\nu_j^2\bar\xi)=-v(\xi)$.
By rotation and $2\cos\theta\,T_j\nu_j=p_j+p_{j+1}$,
\[
  \int_{\partial P}W\,ds
  =\frac{\displaystyle\int_0^1 U((1-t)v_0+tv_1)\cdot\nu_0\,dt}{2\cos\theta}
    \underbrace{\sum_{j=0}^{n-1}|p_{j+1}-p_j|(p_j+p_{j+1})}_{2\int_{\partial P}x\,ds}
      =0 .
\]
Pairing $\xi$ and $\nu_j^2\bar\xi$ in \eqref{sect-matrix} yields
\[
  \begin{aligned}
    & \tfrac12\big(|T_jDU(\xi)T_j^{-1}|^2+|T_jDU(\nu_j^2\bar\xi)T_j^{-1}|^2\big) \\
    & \quad=|DU|^2+(x^2-1)(\partial_{\tau_j}u)^2+(y^2-1)(\partial_{\nu_j}v)^2
      +(xy-1)\big((\partial_{\nu_j}u-\partial_{\tau_j}v)^2
        -2(\partial_{\tau_j}u)(\partial_{\nu_j}v)\big), \\
    & \tfrac12\big(|T_jU(\xi)|^2+|T_jU(\nu_j^2\bar\xi)|^2\big)=(\det T_j)(xu^2+yv^2).
  \end{aligned}
\]
With $x_j,y_j$ as in \eqref{sect-matrix} and $|T_j\tau_j|=\sqrt{y_j\det T_j}$,
integration gives
\begin{equation}
  \label{sect-energy}
  \begin{aligned}
    \frac{n}{2\pi}\int_P|DW|^2\,dX
      &=\sum_{j=0}^{n-1}(\det T_j)
      \big[1+B(x_j^2-1)+C(y_j^2-1)+K(x_jy_j-1)\big], \\
    \frac{s_n}{2\sin\theta}\int_{\partial P}|W|^2\,ds
      &=\sum_{j=0}^{n-1}(\det T_j)^{3/2}\sqrt{y_j}\,[bx_j+(1-b)y_j].
  \end{aligned}
\end{equation}
On each sector write
\[
  \begin{aligned}
    G(x,y) &=1+B(x^2-1)+C(y^2-1)+K(xy-1)>0, & F(x,y) &=bx+(1-b)y>0, \\
    A(r) &=\frac{2\arctan(tr)-2\pi/n}{\sin(2\pi/n)}, &
      m &=2(C-B)-3/2+2b .
  \end{aligned}
\]
\begin{proof}[Proof of Lemma~\ref{sect-scalar}]
  Use $\max_{w>0}\big((G+1/2)w-\sqrt y\,Fw^{3/2}\big)=4(G+1/2)^3/(27yF^2)$.
  Since $\sum_{j=0}^{n-1}\det T_j=n$ and $\sum_{j=0}^{n-1}\phi_j=2\pi$, setting
  $w=\det T_j$ gives
  \[
    \sum_{j=0}^{n-1}(\det T_j)\big(G(x_j,y_j)-|T_j\tau_j|F(x_j,y_j)\big)
    \le \sum_{j=0}^{n-1}\left(\frac{1-\det T_j}{2}
      +m\frac{\phi_j-2\pi/n}{\sin(2\pi/n)}\right)=0 .
  \]
  By \eqref{ray}, \eqref{sect-energy} and $\sigma_1(R_n)=\pi s_n/(n\sin\theta)$,
  \[
    \sigma_1(P)\int_{\partial P}|W|^2\,ds\le\int_P|DW|^2\,dX
    \le\sigma_1(R_n)\int_{\partial P}|W|^2\,ds .
  \]
  Equality requires $(r_j,d_j)=(1,0)$ for $0\le j<n$, hence $p_{j+1}=e^{2\pi i/n}p_j$
  for every $j$.
\end{proof}

\medskip\noindent\emph{Pentagons.}
\begin{theorem}[Equiangular pentagons]
  \label{equiancrefine}
  Every simple equiangular pentagon $P$ satisfies
  \[
    L(P)\sigma_1(P)\le L(R_5)\sigma_1(R_5).
  \]
  Equality holds if and only if $P$ is regular.
\end{theorem}

\begin{theorem}[Local stability of the regular pentagon]
  \label{pentlocaltheorem}
  Let $p_0,\ldots,p_4\in\mathbb C$ be the vertices of $P$ in counterclockwise order.
  Put $v_j=e^{2\pi i j/5}$ for $j\in\mathbb Z$, $p_{j+5}=p_j$ and
  $a=\frac15\sum_{\ell=0}^4p_\ell\overline{v_\ell}$.
  For $a\ne0$ define $d\ge0$ by
  \[
    d^2=\sum_{j=0}^4\left|
      \frac{p_j-\frac15\sum_{\ell=0}^4p_\ell}{a}-v_j\right|^2 .
  \]
  If $a\ne0$ and $d\le1/100$, then
  \begin{equation}
    \label{pen-local-ineq}
    \sqrt{|P|}\,\sigma_1(P)\le\sqrt{|R_5|}\,\sigma_1(R_5)-d^2/100 .
  \end{equation}
\end{theorem}

Akhmetgaliyev, Kao and Osting \cite{AKO} studied fixed-area Steklov optimization
numerically.

\begin{proof}[Proof of Theorem~\ref{equiancrefine}]
  For $z=x_1+ix_2$, $\operatorname{Re}(z)=x_1$, $\operatorname{Im}(z)=x_2$.
  Set $\omega=e^{2\pi i/5}$ and take indices modulo five.
  Then
  \[
    u_z=\tfrac12(\partial_{x_1}u-i\partial_{x_2}u),\qquad
    u_{\bar z}=\tfrac12(\partial_{x_1}u+i\partial_{x_2}u).
  \]

  Put $t=\tan(\pi/5)$ and $k=(\sqrt5-1)/2$.
  Each interior angle is $3\pi/5$.
  Normalize $L(P)=10t$ and the outward normals to $\omega^j$, $0\le j<5$.
  With $h_j=\sup_{z\in P}\operatorname{Re}(\omega^{-j}z)$, side $j$ is
  \[
    z=\omega^j(h_j+iv),\qquad
    \frac{h_j\cos(2\pi/5)-h_{j-1}}{\sin(2\pi/5)}
    \le v\le
    \frac{h_{j+1}-h_j\cos(2\pi/5)}{\sin(2\pi/5)}.
  \]
  Hence $2t\sum_{j=0}^4h_j=L(P)=10t$.
  Translation by $-\frac25\sum_{j=0}^4h_j\omega^j$ gives $\sum_{j=0}^4h_j\omega^j=0$.
  Fourier inversion now yields
  \[
    T=\frac25\sum_{j=0}^4h_j\omega^{-2j},\quad
    h_j=1+\operatorname{Re}(T\omega^{2j}),\quad
    |v+k\operatorname{Im}(T\omega^{2j})|\le t\big(1-(1+k)(h_j-1)\big).
  \]
  For $r=|T|$ and $R_5$ of apothem one,
  \[
    \frac{|P|}{|R_5|}=\frac15\sum_{j=0}^4h_j\big(1-(1+k)(h_j-1)\big)
    =1-\frac{(1+k)r^2}{2}.
  \]
  Since $r\cos(\pi/5)\le\max_{0\le j<5}(h_j-1)<1/(1+k)$, we have $r<0.764$.

  Integration on the sides and $\sum_{j=0}^4\omega^{hj}=0$ for
  $h\in\mathbb Z\setminus5\mathbb Z$ give
  \[
    \frac1{L(P)}\int_{\partial P}z\,ds
    =\frac15\sum_{j=0}^4\omega^j\big(1-(1+k)(h_j-1)\big)(h_j
      -ik\operatorname{Im}(T\omega^{2j}))
    =-\frac{kT^2}{4}.
  \]
  For $u_\alpha(z)=\operatorname{Re}(e^{-i\alpha}(z+kT^2/4))$, integrating
  $2(\operatorname{Re}(w))^2=|w|^2+\operatorname{Re}(w^2)$ on the same intervals gives
  \[
    \frac2{L(P)}\max_{\alpha\in\R}\int_{\partial P}u_\alpha^2\,ds
    =1+\frac{t^2}{3}+(1-2k)r^2-\frac{k^2r^4}{16}
      +\left|\left(2k-\frac{r^2}{2}\right)\overline T-\frac{k^2 T^4}{16}\right|.
  \]
  Since $\int_P|\nabla u_\alpha|^2\,dx=|P|$, \eqref{ray} gives, for
  $2/25\le r\le0.764$,
  \[
    \frac{L(P)\sigma_1(P)}{L(R_5)\sigma_1(R_5)}
    \le \frac{5t}{\pi s_5}
    \frac{1-(1+k)r^2/2}{1+t^2/3+2kr+(1-2k)r^2-r^3/2-k^2r^4/8}<1 .
  \]
  Here $5t/(\pi s_5)<1.232$ by $s_5>47/50$, and concavity bounds the denominator below
  by its smaller endpoint value, greater than $1.27$.
  For $0<r\le2/25$, the reference pentagon has area $|R_5|=5t$ and apothem one.
  Rotate and dilate the eigenfunction in \eqref{eq:regular-eigenfunction} to obtain
  $\varphi(z)=e^{i\pi/5}U_5(e^{-i\pi/5}\pi z/(5t))$.
  Then $\int_{R_5}|\nabla\varphi|^2\,dx=2\pi$ and
  $\nabla\varphi\cdot\overline{\nabla z}=2\varphi_z$.
  Projecting $z$ onto $\mathbb C\varphi$ as in Lemma~\ref{gap} gives
  \[
    \frac1{q_5}\le\frac1{\gamma_5}
    +\left(\frac1{s_5}-\frac1{\gamma_5}\right)
    \frac{\left|\int_{R_5}\varphi_z\,dx\right|^2}{\pi |R_5|}.
  \]
  Since $q_5<0.984<1.48<\gamma_5$, $\int_{R_5}\varphi_z\,dx\ne0$.
  Set $F=|R_5|\varphi/\int_{R_5}\varphi_z\,dx$.
  Since $s_5>0.94$,
  \begin{equation}
    \label{anc-projection}
    \begin{gathered}
      \int_{R_5}(F_z-1)\,dx=0, \\
      \frac1{|R_5|}\int_{R_5}(|F_z-1|^2+|F_{\bar z}|^2)\,dx
        =\frac{\pi |R_5|}{\left|\int_{R_5}\varphi_z\,dx\right|^2}-1
        \le\frac{s_5^{-1}-q_5^{-1}}{q_5^{-1}-\gamma_5^{-1}}<0.14 .
    \end{gathered}
  \end{equation}
  On the $j$th reference triangle, $0\le a\le1$, $|b|\le ta$, define
  \[
    W\bigl(\omega^j(a+ib)\bigr)
      =\omega^j\bigl(ah_j+i(b(1-(1+k)(h_j-1))
        -ak\operatorname{Im}(T\omega^{2j}))\bigr),\qquad V(z)=W(z)-z.
  \]
  The adjacent linear maps agree on their common radius.
  In the normal and tangent frame, with $\xi=V_z$ and $\eta=V_{\bar z}$ in the global
  complex coordinate,
  \[
    DW=
    \begin{pmatrix}
      h_j & 0 \\
      -k\operatorname{Im}(T\omega^{2j}) & 1-(1+k)(h_j-1)
    \end{pmatrix},\quad
    \xi=-kT\omega^{2j}/2,\quad
    \eta=(T\omega^{4j}+(1+k)\overline T)/2 .
  \]
  Thus $\max_{|q|=1}|\xi q+\eta\bar q|=|\xi|+|\eta|<1.62r<1$ on each triangle, so $W$
  is bi-Lipschitz.
  On the boundary,
  \[
    \int_{\partial P}F\circ W^{-1}\,ds
    =\left(\int_{-t}^{t}F(1+ib)\,db\right)
      \sum_{j=0}^4\omega^j\big(1-(1+k)\operatorname{Re}(T\omega^{2j})\big)=0 .
  \]

  Under $z\mapsto z+\varepsilon V(z)$ the arclength factor on side $j$ is
  $1-\varepsilon(1+k)(h_j-1)$.
  For a unit gradient $(\cos\alpha,\sin\alpha)$ the Dirichlet energy factor is
  \[
    \frac{|1+\varepsilon\xi-\varepsilon\eta e^{-2i\alpha}|^2}
      {|1+\varepsilon\xi|^2-\varepsilon^2|\eta|^2}.
  \]
  Its derivative at $\varepsilon=0$ equals $-2\operatorname{Re}(\eta e^{-2i\alpha})$.
  By polarization, the complex bilinear first variation of the Dirichlet energy minus
  $\sigma_1(R_5)$ times the boundary mass is
  \begin{equation}
    \label{anc-derivative}
    \mathcal Q'(u,v)=-4\int_{R_5}(\eta u_zv_z+\bar\eta u_{\bar z}v_{\bar z})\,dx
      +(1+k)\sigma_1(R_5)\sum_{j=0}^4 (h_j-1)\int_{-t}^t(uv)(\omega^j(1+is))\,ds .
  \end{equation}
  By $F_z(\omega z)=F_z(z)$, $F_{\bar z}(\omega z)=\omega^2F_{\bar z}(z)$ and
  $F(\omega z)=\omega F(z)$, one has $\mathcal Q'(F,\bar F)=0$ and
  \[
    \int_{R_5}(\eta F_z^2+\bar\eta F_{\bar z}^2)\,dx
    =\frac{\overline T}{2}\left((1+k)\int_{R_5}F_z^2\,dx
      +\sum_{j=0}^4\omega^{-4j}\int_{S_j}F_{\bar z}^2\,dx\right),
  \]
  where $S_j=\{\omega^j(a+ib):0\le a\le1,\ |b|\le ta\}$.
  By \eqref{anc-projection}, $\int_{R_5}F_z^2\,dx=|R_5|+\int_{R_5}(F_z-1)^2\,dx$,
  hence
  \begin{align*}
    4\left|\int_{R_5}(\eta F_z^2+\bar\eta F_{\bar z}^2)\,dx\right|
      & \ge 2r(1+k)\left(|R_5|-\int_{R_5}\big(|F_z-1|^2+|F_{\bar z}|^2\big)\,dx\right)
        \\
    & \ge1.72(1+k)r|R_5|,
  \end{align*}
  \begin{align*}
    & (1+k)\sigma_1(R_5)\left|\sum_{j=0}^4 (h_j-1)\int_{-t}^t F(\omega^j(1
      +is))^2\,ds\right| \\
    & \qquad\le\frac{(1+k)r}{2}\sigma_1(R_5)\int_{\partial R_5}|F|^2\,ds
      \le 1.14(1+k)r|R_5| .
  \end{align*}
  In the boundary integral only $\bar T\omega^{-2j}/2$ contributes.
  Rotation gives $\int_{R_5}\nabla F\cdot\nabla F\,dx=0$.
  Take $u=\operatorname{Re}(e^{-i\alpha}F)$ with
  $e^{2i\alpha}=-\mathcal Q'(F,F)/|\mathcal Q'(F,F)|$.
  Then
  \[
    \mathcal Q'(u,u)=-|\mathcal Q'(F,F)|/2\le-0.29(1+k)r|R_5|,\qquad
    \int_{R_5}|\nabla u|^2\,dx\le1.14|R_5| .
  \]
  At $\varepsilon=1$, $|\xi|<0.31r$ and $|\eta|<1.31r$ bound the energy remainder by
  \[
    \max_{0\le\alpha\le2\pi}
    \left|\frac{|1-\mu_1 e^{-2i\alpha}|^2}{1-|\mu_1|^2}-1
      +2\operatorname{Re}(\eta e^{-2i\alpha})\right|
    \le\frac{2|\mu_1|^2}{1-|\mu_1|}+\frac{2|\eta\xi|}{|1+\xi|}<5r^2 .
  \]
  The arclength factor is affine in $\varepsilon$.
  Therefore
  \[
    \int_P|\nabla(u\circ W^{-1})|^2\,dx-\sigma_1(R_5)\int_{\partial P}(u\circ
      W^{-1})^2\,ds
    \le |R_5|r\big(5.70r-0.29(1+k)\big)<0 .
  \]
  Thus $L(P)\sigma_1(P)<L(R_5)\sigma_1(R_5)$ for $r>0$, while $r=0$ gives
  $h_0=\cdots=h_4=1$ and $P=R_5$.
\end{proof}

\begin{proof}[Proof of Theorem~\ref{pentlocaltheorem}]
  Use the circumradius-one pentagon, so $|R_5|=\tfrac52\sin(2\pi/5)$ and
  $L(R_5)=10\sin(\pi/5)$.
  If $d=0$, $P$ is similar to $R_5$.
  Otherwise $W_j=d^{-1}((p_j-\frac15\sum_{\ell=0}^4p_\ell)/a-v_j)$ satisfies
  $\sum_{j=0}^4W_j=\sum_{j=0}^4W_j\bar v_j=0$, so
  \[
    W=\sum_{k=2}^4z_kw_k,\qquad w_k=(v_j^k/\sqrt5)_{j=0}^4,\qquad
    z_k=\frac1{\sqrt5}\sum_{j=0}^4W_j\bar v_j^{\,k},
    \qquad \sum_{k=2}^4|z_k|^2=1 .
  \]
  Interpolate each pair $(W_j,W_{j+1})$ by the real linear map
  \begin{equation}
    \label{pen-interpolant}
    B_j(sv_j+tv_{j+1})=sW_j+tW_{j+1}\qquad(s,t\in\R).
  \end{equation}
  Set $V(x)=B_jx$ on the triangle $S_j=\{sv_j+tv_{j+1}:s,t\ge0,\ s+t\le1\}$ and put
  \[
    P_t=\{x+tV(x):x\in R_5\},\qquad Z(x)=-0.45726\,V(x).
  \]
  Here $-0.45726$ is the coefficient of $z$ in the polynomial $f$ of Numerical
  Lemma~\ref{certpentagon}.
  Thus $P_d=(P-\frac15\sum_{\ell=0}^4p_\ell)/a$.
  By \eqref{eq:regular-eigenfunction}, the $L^2(\partial R_5)$ projection of $f$ onto
  the first eigenline in its rotation class is
  \[
    U(x)=\frac{s_5}{L(R_5)}\,U_5\!\left(\frac{2\pi x}{L(R_5)}\right)
      \int_{\partial R_5} f(\xi)\,
      \overline{U_5\!\left(\frac{2\pi\xi}{L(R_5)}\right)}\,ds_\xi
      \qquad(x\in R_5).
  \]
  Indeed $\int_{\partial R_5}|U_5(2\pi\xi/L(R_5))|^2\,ds_\xi=L(R_5)/s_5$.
  In this rotation class Lemma~\ref{gap} places the other eigenvalues above
  $1.48(2\pi/L(R_5))>1.58$.
  The spectral theorem and \eqref{pen-resid} give
  \[
    \min_{\lambda}\,|\lambda-1.017|
    \le \left(\frac{\int_{\partial R_5}|\partial_\nu f-1.017f|^2\,ds}
      {\int_{\partial R_5}|f|^2\,ds}\right)^{1/2}
    <\frac{0.003}{\sqrt{1-10^{-3}}}<0.0031,
  \]
  where $\lambda$ ranges over the Steklov eigenvalues in this class.
  Since $1.017+0.0031<1.58$, the closest eigenvalue is $\sigma_1(R_5)$.
  Expansion of $f-U$ in the remaining eigenfunctions yields
  \begin{align*}
    \left(\int_{\partial R_5}|f-U|^2\,ds\right)^{1/2}
      &<\frac{0.003}{1.58-1.017}<0.0054, \\
    \left(\int_{R_5}|\nabla(f-U)|^2\,dx\right)^{1/2}
      &<\frac{\sqrt{1.58}\,(0.003)}{1.58-1.017}<0.0068 .
  \end{align*}
  By orthogonality on the boundary and in the Dirichlet integral,
  \begin{gather*}
    \int_{\partial R_5}|U|^2\,ds\le\int_{\partial R_5}|f|^2\,ds<(1.001)^2,
      \qquad
      \int_{R_5}|\nabla U|^2\,dx\le\int_{R_5}|\nabla f|^2\,dx<(1.009)^2, \\
    \left|\int_{R_5}|\nabla U|^2\,dx-1.017\right|
      \le\left|\int_{R_5}|\nabla f|^2\,dx-1.017\right|
      +\int_{R_5}|\nabla(f-U)|^2\,dx<0.0003 .
  \end{gather*}

  For $u_t(x+tV(x))=U(x)+tZ(x)$ define
  \[
    D(t)=\sqrt{|P_t|/|R_5|}\int_{P_t}|\nabla u_t|^2\,dx
    -\sigma_1(R_5)\left(\int_{\partial P_t}|u_t|^2\,ds-
    L(P_t)^{-1}\left|\int_{\partial P_t}u_t\,ds\right|^2\right).
  \]
  The symmetries $U(\omega z)=\omega U(z)$, $U(\bar z)=\overline{U(z)}$ hold with
  $\omega=e^{2\pi i/5}$ by the corresponding symmetries of $f$.
  Both $D$ and $|P_t|$ are invariant under
  \[
    (W_j)\mapsto(\omega W_{j-1}),\quad z_k\mapsto\omega^{1-k}z_k,
    \qquad (W_j)\mapsto(\overline{W_{-j}}),\quad z_k\mapsto\bar z_k .
  \]
  Since $\sum_{j=0}^4\omega^{j(1-k)}=0$ for $2\le k\le4$, the rotation average
  annihilates each linear term.
  Thus $D'(0)=\left.\frac{d}{dt}|P_t|\right|_{t=0}=0$ and $D(0)=0$ by \eqref{pde}.

  On $S_j$ and its outer edge put
  \begin{align*}
    G_j &=(\operatorname{tr}B_j)I_2-B_j-B_j^{\mathsf T}, &
      H_j &=B_jB_j^{\mathsf T}-(\det B_j)I_2, \\
    \kappa &=\tfrac15\sum_{j=0}^4\det B_j, \\
    \xi_j &=(W_{j+1}-W_j)/(v_{j+1}-v_j), &
      (\alpha_j,\beta_j) &=(\operatorname{Re}(\xi_j),(\operatorname{Im}(\xi_j))^2/2).
  \end{align*}
  For $J_j=I_2+tB_j$, use $(\det J_j)J_j^{-1}=I_2+t((\operatorname{tr}B_j)I_2-B_j)$ to
  obtain
  \[
    \begin{aligned}
      \det J_j\,J_j^{-1}(J_j^{-1})^{\mathsf T} &=I_2+tG_j+t^2 H_j+O(t^3), \\
      |1+t\xi_j| &=1+t\alpha_j+t^2\beta_j+O(t^3),\qquad
        |P_t|/|R_5|=\tfrac15\sum_{j=0}^4\det J_j=1+\kappa t^2 .
    \end{aligned}
  \]
  Set $(G,H)=(G_j,H_j)$ on $S_j$ and $(\alpha,\beta)=(\alpha_j,\beta_j)$ on its outer
  edge.
  With $\langle a,b\rangle=\sum_{k=1}^2a_k\overline{b_k}$,
  \begin{align}
    \tfrac12D''(0)={} & \int_{R_5}\left(\langle H\nabla U,\nabla U\rangle+
      2\operatorname{Re}\big(\langle G\nabla U,\nabla Z\rangle\big)
        +|\nabla Z|^2\right)\,dx
      +\tfrac\kappa2\int_{R_5}|\nabla U|^2\,dx\notag \\
    & -\sigma_1(R_5)\left(\int_{\partial R_5}
      (\beta|U|^2+2\alpha\operatorname{Re}(U\overline Z)+|Z|^2)\,ds
      -L(R_5)^{-1}\left|\int_{\partial R_5}(\alpha U+Z)\,ds\right|^2\right).
    \label{pen-hess}
  \end{align}
  For the five displacements in \eqref{pen-five}, subtraction of $Q_f$ from
  \eqref{pen-hess} and the bounds \eqref{pen-coeff} give
  \begin{equation}
    \label{pen-er}
    \begin{aligned}
      |D''(0)-Q_f(W)|\le{} &
        7\left(\int_{R_5}|\nabla(U-f)|^2\,dx\right)^{1/2}
        +8\left(\int_{\partial R_5}|U-f|^2\,ds\right)^{1/2} \\
      & +6|\sigma_1(R_5)-1.017|
        +|\kappa|\left|\int_{R_5}|\nabla U|^2\,dx-1.017\right|<0.11 .
    \end{aligned}
  \end{equation}
  Here Cauchy--Schwarz bounds the boundary means, and for a real symmetric matrix $A$
  and $g,h\in\mathbb C^2$,
  \[
    |\langle Ah,h\rangle-\langle Ag,g\rangle|
    \le\max_{|v|=1}|Av|\,(|h|+|g|)|h-g| .
  \]

  The invariant real quadratic forms are spanned by $|z_2|^2$, $|z_3|^2$, $|z_4|^2$
  and $\operatorname{Re}(z_3z_4)$.
  Evaluation at $(w_3\pm w_4)/\sqrt2$ bounds the last coefficient in absolute value by
  $0.75+2(0.11)=0.97$.
  Hence \eqref{pen-five}, \eqref{pen-er} and $2|z_3z_4|\le |z_3|^2+|z_4|^2$ give
  \[
    \begin{aligned}
      D''(0) & \le \sum_{k=2}^4\big(Q_f(w_k)+0.11\big)|z_k|^2
        +0.97\,|\operatorname{Re}(z_3z_4)| \\
      & \le -0.14|z_2|^2-0.18\big(|z_3|^2+|z_4|^2\big)\le-0.14 .
    \end{aligned}
  \]
  By Numerical Lemma~\ref{pent-taylor} and Taylor's formula,
  \[
    D(d)\le(-0.07+6.3d)d^2\le-0.007d^2 .
  \]
  Since $D(d)<0$, Numerical Lemma~\ref{pent-taylor} gives
  \[
    0<\int_{\partial P_d}|u_d|^2\,ds-
      \frac{1}{L(P_d)}\left|\int_{\partial P_d}u_d\,ds\right|^2<1.025 .
  \]
  Applying \eqref{ray} to $\operatorname{Re}(u_d)$ and $\operatorname{Im}(u_d)$ gives
  \[
    \sqrt{|P|}\,\sigma_1(P)-\sqrt{|R_5|}\,\sigma_1(R_5)
      \le\frac{\sqrt{|R_5|}\,D(d)}{1.025}
    \le-\frac{0.007\sqrt{|R_5|}}{1.025}\,d^2<-d^2/100 .
  \]
\end{proof}

\section{The sector comparison and the main theorem}
\label{sec:comparison}

We prove the sector inequality by Bernstein coefficient bounds for $5\le n\le314$ and
uniform estimates for $n\ge315$.
Equality forces all vertex radii to agree and every central angle to equal $2\pi/n$.

\begin{theorem}[Sector comparison]
  \label{sector-comparison}
  Let $n\ge5$, and let $P$ be a convex $n$-gon with
  $\sqrt{|P|}\,\sigma_1(P)\ge\sqrt{|R_n|}\,\sigma_1(R_n)$.
  For the normalized vertices \eqref{eq:normalized-fan}, each sector $0\le j<n$, in
  the coordinates \eqref{sect-matrix}, satisfies
  \begin{equation}
    \label{finite-sector-sign}
    \begin{aligned}
      & 4\big(\tfrac32+B(x^2-1)+C(y^2-1)+K(xy-1)\big)^3 \\
      & \qquad\le27y\big(bx+(1-b)y\big)^2
        \left(\tfrac12+\left(2(C-B)-\tfrac32+2b\right)
        \frac{\phi-2\pi/n}{\sin(2\pi/n)}\right),
    \end{aligned}
  \end{equation}
  with equality precisely when $\phi=2\pi/n$ and $|p_{j+1}|=|p_j|$.
\end{theorem}

\begin{lemma}[Central angles and vertex radii]
  \label{lem:finite-angular-box}
  Under the hypotheses of Theorem~\ref{sector-comparison}, put $h=2\pi/n$ and, for
  each $0\le j<n$,
  \[
    r=\frac{\tan(\tfrac12\arg(p_{j+1}/p_j))}{\tan(\pi/n)},\qquad
    d=\frac{(|p_{j+1}|-|p_j|)^2}{|p_{j+1}|\,|p_j|},\qquad p_n=p_0 .
  \]
  Then
  \begin{equation}
    \label{infinite-sector-geometry}
    \begin{aligned}
      n=5 & : & 0.2<r<3,       & \qquad 0\le d<0.4h^2, \\
      6\le n\le20 & : & 0.155<r<2.20, & \qquad 0\le d<0.23h^2, \\
      21\le n\le314 & : & 0.155<r<1.70, & \qquad 0\le d<0.015h^2, \\
      n\ge315 & : & 0.15<r<1.75, & \qquad 0\le d<3h^4<4\tan^2(\pi/n).
    \end{aligned}
  \end{equation}
\end{lemma}
\begin{proof}
  Normalize $L(P)=2\pi$.
  For $5\le n\le314$ and $c=\pi/(n\tan(\pi/n))$, take the lower bound $s<s_n$ from
  \eqref{eq:certified-regular-bound} and Lemma~\ref{uniformregular},
  \[
    s=
    \begin{cases}
      s_n^{\rm lo} & (5\le n\le19), \\
      1-n^{-3}B(n,16)/(1-k(n,16)) & (20\le n\le314).
    \end{cases}
  \]
  Weinstock's inequality, $\sqrt{|P|}\,\sigma_1(P)\ge\pi s_n/\sqrt{n\tan(\pi/n)}$ and
  polygonal isoperimetry give
  \[
    cs^2=\frac{\pi s^2}{n\tan(\pi/n)}\le\frac{|P|}{\pi}\le\frac{\pi}{n\tan(\pi/n)}=c .
  \]
  For $0\le i<n$, the side length $l_i$ and the normal and tangential coordinates of
  the midpoint are
  \[
    l_i=|p_{i+1}-p_i|,\quad H_i=\frac{\operatorname{Im}(\overline{p_i}p_{i+1})}{l_i},
      \qquad
    m_i=\frac{|p_{i+1}|^2-|p_i|^2}{2l_i},\qquad
    \sum_{i=0}^{n-1}l_i H_i=2|P| .
  \]
  Integration on each side and \eqref{ray} give
  \[
    \frac{2|P|^2}{\pi}+\sum_{i=0}^{n-1}\frac{l_i^3}{12}
    +\sum_{i=0}^{n-1} l_i\big[(H_i-|P|/\pi)^2+m_i^2\big]
    =\int_{\partial P}|x|^2\,ds\le\frac{2|P|}{\sigma_1(P)} .
  \]
  For every $\alpha\ge cs^2/(4\pi)$, Numerical Lemma~\ref{certedge} gives
  \[
    4\pi s^2>3n\tan(\pi/n),\qquad
    \frac{d}{d\alpha}\left(\frac{2\sqrt{n\tan(\pi/n)}}{\pi s}\alpha^{3/2}
      -4\alpha^2\right)
    =\sqrt\alpha\left(\frac{3\sqrt{n\tan(\pi/n)}}{\pi s}-8\sqrt\alpha\right)<0 .
  \]
  Apply this at $\alpha=|P|/(4\pi^2)$ and divide the boundary moment inequality by
  $h^3/12$.
  Then $u_i=l_i/h$ satisfy $\sum_{i=0}^{n-1}u_i=n$ and
  \[
    \sum_{i=0}^{n-1}u_i^3-n+\frac{12}{h^2}\sum_{i=0}^{n-1}u_i[(H_i-|P|/\pi)^2+m_i^2]
    \le\mathcal B_n:=\frac{3n^3}{\pi^2}cs^2(1-cs^2)-n .
  \]
  Fix $0\le i<n$ and write $(u,H,m)=(u_i,H_i,m_i)$.
  Jensen's inequality gives
  \[
    \begin{gathered}
      \sum_{\substack{0\le j<n\\j\ne i}}u_j^3\ge\frac{(n-u)^3}{(n-1)^2},\qquad
        p_n(u):=u^3+\frac{(n-u)^3}{(n-1)^2}-n\le\mathcal B_n, \\
      (H-|P|/\pi)^2+m^2\le e^2:=\frac{h^2}{12u}(\mathcal B_n-p_n(u)).
    \end{gathered}
  \]
  Set $t=\tan(h/2)$, $v=hu/2=uct$ and $E=e/c<s^2$ by Numerical Lemma~\ref{certedge}.
  The endpoint radii satisfy $R_\pm^2=H^2+(m\pm v)^2$.
  Hence
  \[
    H^2-m^2+v^2\le R_-R_+\le H^2+m^2+v^2,\qquad
    \frac{vH}{H^2+m^2}\le\tan(\phi/2)
    =\frac{2Hv}{R_-R_++H^2+m^2-v^2}\le\frac vH .
  \]
  The disk inequality above gives $H\ge|P|/\pi-e>0$ and
  $H^2+m^2\le2|P|H/\pi-|P|^2/\pi^2+e^2\le(|P|/\pi+e)H$.
  Thus
  \begin{equation}
    \label{eq:finite-angular-r}
      \frac{u}{1+E}\le r\le\frac{u}{s^2-E}.
  \end{equation}
  Since $R_+^2-R_-^2=4mv$ and $R_-^2R_+^2=(H^2+m^2+v^2)^2-4v^2m^2$,
  \[
    d=2\left(\sqrt{1+\frac{4v^2m^2}{R_-^2 R_+^2}}-1\right).
  \]
  At fixed $|m|$, the ratio $|m|/(H^2+m^2+v^2)$ decreases in $H>0$.
  On the bounding circle its maximum is
  \[
    \max_{-1\le y\le1}\frac{e\sqrt{1-y^2}}{|P|^2/\pi^2+e^2+v^2+2|P|ey/\pi}
    =\frac{e}{\sqrt{(|P|^2/\pi^2-e^2+v^2)^2+4e^2v^2}} .
  \]
  Therefore $2v|m|/(R_-R_+)\le2ve/(|P|^2/\pi^2-e^2+v^2)$.
  With $|P|/\pi\ge cs^2$,
  \begin{equation}
    \label{eq:finite-angular-d}
    d\le2\left(\sqrt{1+\left(\frac{2utE}{s^4-E^2+u^2t^2}\right)^2}-1\right).
  \end{equation}
  For $5\le n\le314$ the bounds follow from Numerical Lemma~\ref{certedge}.

  \medskip\noindent\emph{The range $n\ge315$.}
  Put $\theta=\pi/n$, $h=2\theta$, $p=\tan^2\theta$ and $c=\theta\cot\theta$.
  By $\theta<\tan\theta<\theta/\cos\theta$ and Lemma~\ref{uniformregular},
  \begin{equation}
    \label{infinite-regular}
    h<0.02,\qquad 9.869<n^2p<9.871,\qquad s_n>s:=1-2.435/n^3 .
  \end{equation}
  The alternating series at $0<\theta<0.01$ give
  \[
    c\ge\frac{1-\theta^2/2+\theta^4/24-\theta^6/720}
      {1-\theta^2/6+\theta^4/120}
    \ge1-\theta^2/3-\theta^4/40,
    \qquad 0<1-cs^2\le\theta^2/3+0.158\theta^3+0.025\theta^4 .
  \]
  Since $cs^2>0.99996$ and $4cs^2>3$, the same midpoint calculation gives, for
  $u=|p_{i+1}-p_i|/h$, $0\le i<n$,
  \[
    (u-1)^2(u+2)\le p_n(u)\le\mathcal B_n
    =\frac{3n^3}{\pi^2}cs^2(1-cs^2)-n\le3\pi(0.158+0.025\theta)<1.50 .
  \]
  Thus $0.16<u<1.7$ and $e=h\sqrt{(\mathcal B_n-p_n(u))/(12u)}<0.885h$.
  With $e/c<0.01771$, \eqref{eq:finite-angular-r}--\eqref{eq:finite-angular-d} yield
  \[
    0.15<r<1.75,\qquad
    0\le d\le\frac{(hue)^2}{(c^2s^4-e^2)^2}<3h^4<4p .
  \]
\end{proof}

\begin{lemma}
  \label{gen-sector-projection}
  Let $n\ge5$, $S_n\le s_n$ and $h=\pi/n$.
  With $q_n,\gamma_n$ given by \eqref{linearray} and Lemma~\ref{gap}, set
  \[
    \eta=\frac{1-S_n/q_n}{1-S_n/\gamma_n},\qquad
    \epsilon=\sqrt{\frac{q_n-S_n}{\gamma_n-S_n}}.
  \]
  If $0<\eta<1$, then
  \[
    B,C\ge0,\qquad B+C\le\frac{\eta}{1-\eta},\qquad 0\le K+B+C
      \le\frac{2\eta}{1-\eta} .
  \]
  If $\epsilon^2<\tan^2h/(3+\tan^2h)$, then
  \[
    \frac{\big(\tan h\sqrt{1-\epsilon^2}-\sqrt3\,\epsilon\big)^2}{3+\tan^2 h}
    \le 1-b\le
    \frac{\big(\tan h\sqrt{1-\epsilon^2}+\sqrt3\,\epsilon\big)^2}{3+\tan^2 h} .
  \]
\end{lemma}
\begin{proof}
  At circumradius one the boundary projection of the coordinate function $x\mapsto x$
  onto the $+1$ first eigenline is
  \[
    F(x)=\frac{s_n U(x)}{2n\sin h}\int_{\partial R_n}\xi\,\overline{U(\xi)}\,ds_\xi .
  \]
  Reflection makes the coefficient of $U$ real.
  Testing the eigenfunction identity by $x-F(x)$ gives
  \[
    \int_{R_n}\nabla(x-F(x))\cdot\overline{\nabla F(x)}\,dx
    =\sigma_1(R_n)\int_{\partial R_n}(x-F(x))\overline{F(x)}\,ds=0 .
  \]
  Both $F$ and $x-F(x)$ have rotation class $+1$.
  Write $|I_2-DF|^2=|1-\partial_{x_1}F|^2+|i-\partial_{x_2}F|^2$ in the complex
  coordinate $x=x_1+ix_2$.
  Then \eqref{linearray} and Lemma~\ref{gap} give
  \[
    \begin{aligned}
      \frac{2|R_n|}{q_n}
        &=\frac{2\pi}{L(R_n)}\int_{\partial R_n}|x|^2\,ds
        \le\frac{2|R_n|}{s_n}
        -\left(\frac1{s_n}-\frac1{\gamma_n}\right)\int_{R_n}|I_2-DF|^2\,dx, \\
      \frac{\int_{R_n}|I_2-DF|^2\,dx}{2|R_n|}
        & \le\frac{1-s_n/q_n}{1-s_n/\gamma_n}\le\eta,\qquad
        \frac{\int_{\partial R_n}|x-F(x)|^2\,ds}{\int_{\partial R_n}|x|^2\,ds}
        \le\frac{q_n}{\gamma_n}\eta=\epsilon^2 .
    \end{aligned}
  \]
  Here $S_n\le s_n\le q_n<\gamma_n$.
  Rotation gives both orthogonalities and ratios on each sector.
  Since $\eta<1$, $F$ is a nonzero real multiple of $U$.
  For $u=F\cdot\nu_j$, $v=F\cdot\tau_j$, the entries of $I_2-DF$ give
  \[
    \begin{aligned}
      (\partial_{\tau_j}u)^2+(\partial_{\nu_j}v)^2 & \le|I_2-DF|^2, \\
      (\partial_{\nu_j}u-\partial_{\tau_j}v)^2+(\partial_{\tau_j}u
        -\partial_{\nu_j}v)^2 & \le2|I_2-DF|^2,\qquad
        \frac{\int_{S_j}|I_2-DF|^2\,dx}{\int_{S_j}|DF|^2\,dx}\le\frac{\eta}{1-\eta}.
    \end{aligned}
  \]
  Integration gives the bounds on $B+C$ and $K+B+C$.
  On $\Gamma_j$ the normal mass fractions are
  \[
    b=\frac{\int_{\Gamma_j}(F\cdot\nu_j)^2\,ds}{\int_{\Gamma_j}|F|^2\,ds},
    \qquad \frac{\int_{\Gamma_j}(x\cdot\nu_j)^2\,ds}{\int_{\Gamma_j}|x|^2\,ds}
      =\frac{3}{3+\tan^2 h}.
  \]
  By orthogonality and Cauchy--Schwarz in the two components,
  \[
    \sqrt{1-\epsilon^2}\le
    \sqrt{1-\frac{\int_{\Gamma_j}|x-F(x)|^2\,ds}{\int_{\Gamma_j}|x|^2\,ds}}
    =\frac{\int_{\Gamma_j}F\cdot x\,ds}
      {\sqrt{\int_{\Gamma_j}|F|^2\,ds\int_{\Gamma_j}|x|^2\,ds}}
    \le\frac{\sqrt{3b}+\tan h\sqrt{1-b}}{\sqrt{3+\tan^2h}} .
  \]
  Hence
  \[
    \left|\arcsin\sqrt b-\arctan\frac{\sqrt3}{\tan h}\right|\le\arcsin\epsilon,
    \qquad 0\le\arctan\frac{\sqrt3}{\tan h}\pm\arcsin\epsilon\le\frac\pi2,
  \]
  where $h\le\pi/5$ and $\epsilon^2<\tan^2h/(3+\tan^2h)$.
  Taking squared cosines gives the bounds.
\end{proof}

For $5\le n\le20$ we compare the moments of the first eigenfunction with those of a
harmonic polynomial.
Orient the circumradius-one polygon by
\[
  z(t)=\cos h+it\sin h,\quad E=\{z(t):|t|\le1\},\quad
  S=\{rz(t):0\le r\le1,\ |t|\le1\},\qquad h=\pi/n .
\]
For a nonzero harmonic polynomial $f$ with $f(e^{2ih}z)=e^{2ih}f(z)$ and
$f(\bar z)=\overline{f(z)}$, the boundary projection is
\[
  U(z)=\frac{s_n}{2n\sin h}\,U_n\!\left(\frac{\pi e^{ih}z}{n\sin h}\right)
    \int_{\partial R_n}f(\xi)\,
    \overline{U_n\!\left(\frac{\pi e^{ih}\xi}{n\sin h}\right)}\,ds_\xi .
\]
The line is spanned by $e^{-ih}U_n(\pi e^{ih}z/(n\sin h))$, which commutes with
conjugation.
Thus $U(\bar z)=\overline{U(z)}$.
The outward normal to $E$ is $(1,0)$.
Put
\[
  g=\frac{(n-1)\cos h}{(\sqrt{1+\tan^2 h/4}+\tan h/2)^2}.
\]
For $0<\mu<g$, Lemma~\ref{gap}, spectral expansion and rotation give
\begin{equation}
  \label{eq:sector-trial-residual}
  \begin{aligned}
    \frac{\int_E|f-U|^2\,ds}{\int_E|f|^2\,ds}
      & \le e_b^2:=
      \frac{\int_E|\partial_x f-\mu f|^2\,ds}
      {(g-\mu)^2\int_E|f|^2\,ds}, \\
    \frac{\int_S\big(|\partial_x(f-U)|^2+|\partial_y(f-U)|^2\big)\,dx\,dy}
      {\int_S(|\partial_x f|^2+|\partial_y f|^2)\,dx\,dy}
      & \le e_g^2:=
      \frac{g\int_E|\partial_x f-\mu f|^2\,ds}
      {(g-\mu)^2\int_S(|\partial_x f|^2+|\partial_y f|^2)\,dx\,dy}.
  \end{aligned}
\end{equation}
Here $e_b,e_g\ge0$ and $t\ge g>\mu>0$ gives
\[
  \frac{d}{dt}\frac{1}{t-\mu}=-\frac{1}{(t-\mu)^2}<0,\qquad
  \frac{d}{dt}\frac{\sqrt t}{t-\mu}=-\frac{t+\mu}{2\sqrt t(t-\mu)^2}<0 .
\]
The normalized moments on $E$ and $S$ are
\[
  \begin{aligned}
    b_f &=\frac{\int_E (\operatorname{Re}(f))^2\,ds}{\int_E|f|^2\,ds}, \\
    B_f &=\frac{\int_S (\operatorname{Re}(\partial_y f))^2\,dx\,dy}
      {\int_S |\nabla f|^2\,dx\,dy},\qquad
      C_f=\frac{\int_S (\operatorname{Im}(\partial_x f))^2\,dx\,dy}
        {\int_S |\nabla f|^2\,dx\,dy}, \\
    K_f &=\frac{\displaystyle\int_S\!
      \big[(\operatorname{Re}(\partial_x f)-\operatorname{Im}(\partial_y f))^2
      -2\operatorname{Re}(\partial_y f)\operatorname{Im}(\partial_x f)\big]\,dx\,dy}
      {\displaystyle\int_S |\nabla f|^2\,dx\,dy}.
  \end{aligned}
\]
If $e_b,e_g<1$, Pythagoras and the triangle inequality give
\[
  \begin{aligned}
    \max\{\sqrt{1-b_f}-e_b,0\}^2 & \le 1-b\le
      \frac{(\sqrt{1-b_f}+e_b)^2}{1-e_b^2}, \\
    \max\{\sqrt{B_f}-e_g,0\}^2 & \le B\le
      \frac{(\sqrt{B_f}+e_g)^2}{1-e_g^2}, \\
    \max\{\sqrt{C_f}-e_g,0\}^2 & \le C\le
      \frac{(\sqrt{C_f}+e_g)^2}{1-e_g^2}.
  \end{aligned}
\]
For $K$, expansion of \eqref{sect-moments} gives
\[
  \begin{aligned}
    & \left|
      K\frac{\int_S(|\partial_x U|^2+|\partial_y U|^2)\,dx\,dy}
      {\int_S(|\partial_x f|^2+|\partial_y f|^2)\,dx\,dy}-K_f\right| \\
    & \quad\le 2e_g\left[
      \left(\frac{2\int_S(\partial_x\operatorname{Re}(f)
        -\partial_y\operatorname{Im}(f))^2\,dx\,dy}
      {\int_S(|\partial_x f|^2+|\partial_y f|^2)\,dx\,dy}\right)^{1/2}
      +\sqrt{B_f}+\sqrt{C_f}\right]+4e_g^2.
  \end{aligned}
\]
Energy orthogonality gives
\[
  (1-e_g^2)\!\int_S(|\partial_x f|^2+|\partial_y f|^2)\,dx\,dy
  \le\!\int_S(|\partial_x U|^2+|\partial_y U|^2)\,dx\,dy
  \le\!\int_S(|\partial_x f|^2+|\partial_y f|^2)\,dx\,dy .
\]

\begin{proof}[Proof of Theorem~\ref{sector-comparison}]
  \emph{The range $5\le n\le314$.} Let $h=2\pi/n$, $p=\tan^2(\pi/n)$.
  Lemma~\ref{lem:finite-angular-box} places every $(r,d)$ in
  $[r_{\min},r_{\max}]\times[0,d_{\max}]$, where
  \[
    (r_{\min},r_{\max},d_{\max}/h^2)=
    \begin{cases}
      (1/5,3,2/5), & n=5, \\
      (1/8,13/5,201/100), & 6\le n\le20, \\
      (7/50,9/5,2), & 21\le n\le314 .
    \end{cases}
  \]
  For $5\le n\le20$ use the $16$ vertices of \eqref{eq:finite-scaled-moments}.
  For $21\le n\le314$, Lemma~\ref{gen-sector-projection} and Numerical
  Lemma~\ref{sector-spectral-intervals}, with $S_n=1-3.044/n^3$, give the $12$
  vertices
  \begin{equation}
    \label{eq:coefficient-vertices}
    \begin{gathered}
      b\in\left\{1-\frac{\big(\tan(\pi/n)\sqrt{1-\epsilon^2}
        \mathbin{\pm}\sqrt3\,\epsilon\big)^2}{3+\tan^2(\pi/n)}\right\}, \\
      (B,C)\in\{(0,0),(\frac{\eta}{1-\eta},0),(0,\frac{\eta}{1-\eta})\},\qquad
        K\in\{-B-C,\,\frac{2\eta}{1-\eta}-B-C\}.
    \end{gathered}
  \end{equation}
  Fix a coefficient vertex and $d\in\{0,d_{\max}\}$.
  We use
  \[
    A(r)=\frac{1+p}{\sqrt p}\arctan\!\left(\frac{\sqrt p(r-1)}{1+pr}\right)
    \ge\frac{V(r)}{T(r)},\qquad T(r)>0,
  \]
  where, for $r\le1$,
  \[
    T(r)=15(1+pr)^5,\quad
    V(r)=(1+p)(r-1)\big[15(1+pr)^4-5p(r-1)^2(1+pr)^2+3p^2(r-1)^4\big],
  \]
  and, for $r\ge1$,
  \[
    T(r)=3(1+pr)^2+p(r-1)^2,\qquad V(r)=3(1+p)(r-1)(1+pr) .
  \]
  For $z\le0$, $\arctan z\ge z-z^3/3+z^5/5$.
  For $z\ge0$, $\arctan z=z\int_0^1(1+z^2s^2)^{-1}\,ds\ge z/(1+z^2/3)$.

  On either side of $r=1$ the expression
  \begin{equation}
    \label{eq:sector-polynomial-defect}
    P(r)=Tr^6\big[27(1/2+mV/T)yF^2-4(G+1/2)^3\big]
  \end{equation}
  is a polynomial by \eqref{sect-matrix} and \eqref{sect-energy}.

  At $d=0$, $y(1)=F(1)=G(1)=1$ and $(4(G+1/2)^3/(27yF^2))'(1)=2(C-B)-3/2+2b=m$.
  Both branches satisfy $(V/T)(1)=0$, $(V/T)'(1)=1$, so $P(1)=P'(1)=0$.

  On the left take $(\delta,l)=(-1,1-r_{\min})$, and on the right
  $(\delta,l)=(1,r_{\max}-1)$.
  Put $j=2$ at $d=0$ and $j=0$ at $d=d_{\max}$.
  With $D=17-j$ on the left and $D=14-j$ on the right, the Bernstein coefficients of
  $P(1+\delta w)/w^j$ on $0\le w\le l$ are
  \begin{equation}
    \label{eq:explicit-bernstein-coefficients}
    \beta_k=\sum_{i=0}^k
    \frac{\binom{k}{i}}{\binom{D}{i}}l^i
    \frac{\delta^{i+j}P^{(i+j)}(1)}{(i+j)!},\qquad 0\le k\le D .
  \end{equation}

  Numerical Lemma~\ref{certbern} gives $m>0$, $1/2+m(V/T)(r_{\min})>0$ and $\beta_k>0$
  for $0\le k\le D$.
  For $0\le w\le l$, the Bernstein expansion \cite{Farouki} gives
  \[
    \frac{P(1+\delta w)}{w^j}
    =\sum_{k=0}^D\beta_k\binom Dk (w/l)^k(1-w/l)^{D-k}
    \ge\min_{0\le k\le D}\beta_k>0 ,
  \]
  with the quotient extended polynomially to $w=0$.
  Hence $P(r)=0$ precisely when $(r,d)=(1,0)$.

  Since the left minorant increases and the right minorant is nonnegative,
  $1/2+mA(r)>0$.
  From $m>0$, $A(r)\ge V/T$ and \eqref{eq:sector-polynomial-defect}, we obtain
  \eqref{finite-sector-sign} at both endpoints, strictly unless $(r,d)=(1,0)$.
  At every coefficient vertex $0<b<1$.

  For positive affine functions $X,Y$ of one variable,
  \[
    (X^{1/3}Y^{2/3})''=-\tfrac29 X^{1/3}Y^{2/3}(X'/X-Y'/Y)^2\le0 .
  \]
  At either $d$-endpoint take $X=1/2+mA(r)$ and $Y=F$ on each line segment of moments.
  Both are positive on the hull by affinity.
  Since $G$ is affine in $(b,B,C,K)$, the cube-root inequality
  \[
    \left(\tfrac{27}4\big(1/2+mA(r)\big)yF^2\right)^{1/3}-G-1/2\ge0
  \]
  extends to the entire hull.
  Now fix the moments of $U$ and take $X=y$, $Y=F$ as affine functions of $d$.
  By \eqref{sect-matrix} and Cauchy--Schwarz in \eqref{sect-moments},
  \[
    \partial_d y=p^{-1}\partial_d x,\qquad K\ge-2\sqrt{BC},\qquad
    \tfrac12\partial_d^2G=(\partial_d x)^2(B+K/p+C/p^2)
    \ge(\partial_d x)^2(\sqrt B-\sqrt C/p)^2 .
  \]
  Concavity on $[0,d_{\max}]$ yields \eqref{finite-sector-sign}, strictly for $d>0$ or
  $r\ne1$.

  \medskip\noindent\emph{The range $n\ge315$.}
  Put $\theta=\pi/n$ and $p=\tan^2\theta$ as in \eqref{infinite-regular}.
  By Lemma~\ref{gap}, $\tan\theta<0.010001$, $\sin\theta/\theta\ge1-\theta^2/6$ and
  $\cos\theta\ge1-\theta^2/2$ give $\gamma_n>0.98n$.
  From $1-2.435/n^3<s_n\le q_n\le1$ and Lemma~\ref{gen-sector-projection} with
  $S_n=1-2.435/n^3$,
  \[
    \eta\le\frac{2.435/n^3}{1-1/(0.98n)}<\frac{2.45}{n^3},\qquad
    \frac{\eta}{1-\eta}<\frac{2.46}{n^3},\qquad \epsilon<\frac{1.59}{n^2} .
  \]
  Since $1.813<n\tan\theta/\sqrt{3+\tan^2\theta}<1.815$,
  \[
    \epsilon^2<\frac{\tan^2\theta}{3+\tan^2\theta},\qquad
    1.813\sqrt{1-(1.59/n^2)^2}-1.59/n\le n\sqrt{1-b}\le1.815+1.59/n .
  \]
  Hence
  \begin{equation}
    \label{infinite-sector-spectral}
    3.24<n^2(1-b)<3.34,\qquad
    0\le n^3B,n^3C<2.46,\qquad n^3K<5 .
  \end{equation}

  Write $D(r,d)=4(G+1/2)^3/(27yF^2)-1/2-mA(r)$.
  By \eqref{infinite-sector-spectral}, $0.4999<m=1/2+2(C-B)-2(1-b)<1/2$.
  Integrating $1-A'(r)=p(r^2-1)/(1+pr^2)$ from $1$ to $r$ gives
  \begin{equation}
    \label{infinite-slope}
    0\le r-1-A(r)
    =p(r-1)^2\int_0^1\frac{\tau(2+\tau(r-1))}{1+p(1+\tau(r-1))^2}\,d\tau
    \le \tfrac13p(r+2)(r-1)^2 .
  \end{equation}
  Thus $1/2+mA(r)>r/2-0.001>0$.
  By the same concavity in $d$, it suffices to bound $D(r,4p)$ and $D(r,0)$.

  At $d=4p$, \eqref{sect-matrix} and \eqref{infinite-sector-spectral} give
  \[
    x=(1+p+p^2r^2)/r<7,\quad r+1/r\le y=r+1/r+pr<9,\quad
    rF\ge b>0.99,\quad G<1.001 .
  \]
  Consequently
  \[
    D(r,4p)<\frac{0.52r^3}{1+r^2}-\frac r2+0.001
    =0.001-\frac{r(1-0.04r^2)}{2(1+r^2)}<0
    \qquad(0.15<r<1.75).
  \]

  At $d=0$, $x=1/r$, $y=r$, $rF=1+(1-b)(r^2-1)$ and
  \[
    \frac{2G+1}{3}-rF=(r^2-1)\left(\tfrac23(C-B/r^2)-(1-b)\right).
  \]
  Expansion about $rF$ gives
  \begin{equation}
    \label{infinite-divided}
    \begin{aligned}
      D(r,0)={} & (r-1)^2\big((C-1+b)(r+2)+B/r\big)+m(r-1-A(r)) \\
      & +\frac{\big((2G+1)/3-rF\big)^2\big((2G+1)/3+2rF\big)}{2rF^2}.
    \end{aligned}
  \end{equation}
  By \eqref{infinite-sector-geometry} and \eqref{infinite-sector-spectral},
  \[
    rF>0.99,\qquad \left|\tfrac{2G+1}{3}-rF\right|\le4|r^2-1|/n^2,
    \qquad (2G+1)/(3rF)+2<4 .
  \]
  The last term in \eqref{infinite-divided} is at most
  $32(r+1)^2(r-1)^2/(n^4F)<0.01(r-1)^2/n^2$ for $r\ne1$.
  Substitution of $n^2(C-1+b)<-3.2$, $n^2B/r<0.06$ and \eqref{infinite-slope} gives
  \[
    \frac{n^2D(r,0)}{(r-1)^2}
    <\left(\frac{9.871}{6}-3.2\right)(r+2)+0.07<0
    \qquad(r\ne1).
  \]
  Together with $D(1,0)=0$ this proves \eqref{finite-sector-sign}, with equality
  precisely at $(r,d)=(1,0)$.
\end{proof}

\begin{proof}[Proof of Theorem~\ref{main}]
  Assume $\sqrt{|P|}\,\sigma_1(P)\ge\sqrt{|R_N|}\,\sigma_1(R_N)$.
  By \eqref{eq:convex-reduction} and \eqref{eq:normalized-fan}, $P$ is a convex
  $N$-gon normalized by
  \[
    |P|=|R_N|,\qquad \int_{\partial P}x\,ds=0,\qquad \sigma_1(P)\ge\sigma_1(R_N).
  \]
  For $W(T_j\xi)=T_jU(\xi)$, Lemma~\ref{sect-scalar}, Theorem~\ref{sector-comparison}
  and \eqref{ray} give $\int_{\partial P}W\,ds=0$, $\int_{\partial P}|W|^2\,ds>0$ and
  \[
    \sigma_1(R_N)\le\sigma_1(P)
    \le\frac{\displaystyle\int_P(|\partial_x W|^2+|\partial_y W|^2)\,dx\,dy}
      {\displaystyle\int_{\partial P}|W|^2\,ds}
    \le\sigma_1(R_N).
  \]
  Equality in \eqref{sect-sign} on each sector gives
  \[
    \phi_j=2\pi/N,\quad |p_{j+1}|=|p_j|
    \quad\Longrightarrow\quad p_{j+1}=e^{2\pi i/N}p_j\qquad(0\le j<N).
  \]
  In the original coordinates, with $a_N=a_0$,
  \begin{equation}
    \label{eq:regular-vertices}
    a_j=a_0e^{2\pi i j/N}+(1-e^{2\pi i j/N})
    \frac{\sum_{i=0}^{N-1}|a_{i+1}-a_i|(a_i+a_{i+1})}
      {2\sum_{i=0}^{N-1}|a_{i+1}-a_i|},\qquad 0\le j<N .
  \end{equation}
  Conversely, \eqref{eq:regular-vertices} gives equality by similarity invariance.
\end{proof}

\appendix
\addtolength{\abovedisplayskip}{-1pt}
\addtolength{\belowdisplayskip}{-1pt}

\section{Finite numerical estimates}
\label{app:numerical}

Rational interval estimates bound the sector moments, prove positivity of the
Bernstein coefficients, and control perturbations of the regular pentagon.

\begin{numericallemma}[Positivity of the Bernstein coefficients]
  \label{certbern}
  Let $5\le n\le314$.
  For each coefficient vertex in \eqref{eq:finite-scaled-moments} or
  \eqref{eq:coefficient-vertices}, each $d\in\{0,d_{\max}\}$ and either branch of
  \eqref{eq:explicit-bernstein-coefficients}, one has
  \[
    m>0.04,\qquad \frac12+m\frac{V(r_{\min})}{T(r_{\min})}>0.07,\qquad
    \beta_k>7\cdot10^{-7}\quad(0\le k\le D),
  \]
  where $V(r_{\min})/T(r_{\min})$ uses the branch $r\le1$.
\end{numericallemma}

\begin{numericallemma}[Spectral moments of a regular sector]
  \label{sector-spectral-intervals}
  The sector moments \eqref{sect-moments} satisfy, for $5\le n\le20$,
  \begin{equation}
    \label{eq:finite-scaled-moments}
    \tfrac72\le n^2(1-b)\le\tfrac{27}5,\quad
    0\le n^3B\le\tfrac85,\quad 0\le n^3C\le\tfrac{13}5,\quad
    -2\le n^3K\le\tfrac52 .
  \end{equation}
  For $21\le n\le314$, the choice $S_n=1-3.044/n^3$ in
  Lemma~\ref{gen-sector-projection} satisfies $0<\eta<3.24/n^3$ and
  $\epsilon^2<\tan^2(\pi/n)/(3+\tan^2(\pi/n))$.
\end{numericallemma}

All decimals are rational.
For $x\in[a,b]$, $y\in[c,d]$, $Q\in\mathbb N$ and $\circ\in\{+,-,\times,/\}$, outward
rounding \cite{MKC} gives
\[
  \frac1Q\min_{\substack{v\in\{a,b\}\\w\in\{c,d\}}}\lfloor Q(v\circ w)\rfloor
  \le x\circ y\le
  \frac1Q\max_{\substack{v\in\{a,b\}\\w\in\{c,d\}}}\lceil Q(v\circ w)\rceil .
\]
Division requires $cd>0$.
For $a\ge0$ the square-root bounds are
$Q^{-1}\lfloor Q\sqrt a\rfloor\le\sqrt x\le Q^{-1}(1+\lfloor Q\sqrt b\rfloor)$.

$\pi=16\arctan(1/5)-4\arctan(1/239)$ and the alternating arctangent series give
$\pi_-<\pi<\pi_+$ for
\[
  \pi_-=4\sum_{i=0}^{51}\frac{(-1)^i}{2i+1}
    \left(\frac4{5^{2i+1}}-\frac1{239^{2i+1}}\right),\qquad
  \pi_+=\pi_-+\frac{16}{105\cdot5^{105}} .
\]
The alternating series for sine and cosine give
\begin{equation}
  \label{eq:rational-trig-endpoints}
  \frac{\displaystyle\sum_{j=0}^{11}\frac{(-1)^j(\pi_-/n)^{2j+1}}{(2j+1)!}}
    {\displaystyle\sum_{j=0}^{10}\frac{(-1)^j(\pi_-/n)^{2j}}{(2j)!}}
  =:t_-(n)<\tan(\pi/n)<t_+(n):=
  \frac{\displaystyle\sum_{j=0}^{10}\frac{(-1)^j(\pi_+/n)^{2j+1}}{(2j+1)!}}
    {\displaystyle\sum_{j=0}^{11}\frac{(-1)^j(\pi_+/n)^{2j}}{(2j)!}} .
\end{equation}
\begin{numericallemma}[Schur bounds for the regular polygon]
  \label{finite}
  Normalize each regular polygon by $L(R_n)=2\pi$.
  Then
  \[
    \sigma_1(R_5)>47/50,\qquad \sigma_1(R_6)>243/250,\qquad
    \sigma_1(R_n)^2>\frac{n\tan(\pi/n)}{(n-1)\tan(\pi/(n-1))}
    \quad(7\le n\le19).
  \]
\end{numericallemma}
\begin{proof}
  The Fourier recurrence gives the integer bounds
  \begin{equation}
    \label{eq:schur-recursion}
    \begin{gathered}
      V_0=10^{12},\qquad
        10^{12}\rho_{n,j+1}\le V_{j+1}:=\left\lceil\frac{(nj+1)V_j}
          {nj+n-1}\right\rceil
        \quad(0\le j<512), \\
      \frac{10^{12}\rho_{n,|i|}}{|1+ni|}\le C_i:=
        \left\lceil\frac{V_{|i|}}{|1+ni|}\right\rceil\qquad(1\le|i|\le256).
    \end{gathered}
  \end{equation}
  Substitution in $k,B,H$ and bisection in \eqref{eq:certified-regular-bound} give
  \[
    s_5^{\rm lo}>0.945,\qquad s_6^{\rm lo}>0.975,\qquad
    \min_{7\le n\le19}\!
    \left((s_n^{\rm lo})^2-\frac{n\,t_+(n)}{(n-1)t_-(n-1)}\right)>10^{-4},
  \]
  where \eqref{eq:rational-trig-endpoints} gives $t_-(n)<\tan(\pi/n)<t_+(n)$.
\end{proof}

\begin{numericallemma}[Bounds for a single edge]
  \label{certedge}
  For $5\le n\le314$ take $s$ from the proof of Lemma~\ref{lem:finite-angular-box} and
  put $h=2\pi/n$, $t=\tan(\pi/n)$.
  Then $4\pi s^2>3nt$ and
  \[
    1.470<\mathcal B_n=\frac{3n^2s^2}{t}\left(\frac1\pi-\frac{s^2}{nt}\right)-n<1.787.
  \]
  For $u>0$ such that $p_n(u)=u^3+(n-u)^3/(n-1)^2-n\le\mathcal B_n$, put
  \[
    E=2t\sqrt{\frac{\mathcal B_n-p_n(u)}{12u}}.
  \]
  Then $E<s^2$ and the triple
  \[
    \left(\frac{u}{1+E},\ \frac{u}{s^2-E},
      \ \frac{2}{h^2}\left[\sqrt{1+\left(\frac{2utE}{s^4-E^2+u^2t^2}\right)^2}
      -1\right]\right)
  \]
  belongs, respectively, to
  \[
    \begin{cases}
      [0.231,\infty)\times(0,2.857]\times[0,0.385], & n=5, \\
      [0.163,\infty)\times(0,2.163]\times[0,0.218], & 6\le n\le20, \\
      [0.164,\infty)\times(0,1.671]\times[0,0.014], & 21\le n\le314 .
    \end{cases}
  \]
\end{numericallemma}
\begin{proof}
  Write $t_\pm=t_\pm(n)$ and $y_\pm=s^2/(nt_\mp)$.
  Rational substitution gives
  \[
    4\pi_-s^2>3nt_+,\qquad
    1.470<3n^3y_+(1/\pi_+-y_+)-n
    \le \mathcal B_n\le 3n^3y_-(1/\pi_--y_-)-n<1.787 .
  \]
  Indeed $1/\pi-2y_-<0$.

  Since $p_n(1)=0$, $p_n(2)>2$ and $p_n'(u)=3n(u-1)((n-2)u+n)/(n-1)^2$, it suffices to
  examine $[a,b]=[i/2000,(i+1)/2000]$ for $0\le i<4000$.
  Since $\min_{a\le u\le b}p_n(u)=p_n(\min\{b,\max\{1,a\}\})$, take $Q=10^{25}$ and
  \[
    \Delta_{ni}=\frac{\left\lceil Q\big(3n^3y_-(1/\pi_--y_-)-n\big)\right\rceil
      -\lfloor Qp_n(\min\{b,\max\{1,a\}\})\rfloor}{Q}.
  \]
  For $\Delta_{ni}\ge0$, one has $a>0$ and \eqref{eq:finite-angular-d} gives
  \[
    E\le\widehat E=\frac{2\lceil Q t_+\rceil}{Q}\sqrt{\frac{\Delta_{ni}}{12a}},
    \qquad
    \frac{2utE}{s^4-E^2+u^2t^2}\le\widehat Z=
    \frac{2b\,\lceil Qt_+\rceil\widehat E/Q}
    {\lfloor Qs^4\rfloor/Q-\widehat E^2+(a\lfloor Qt_-\rfloor/Q)^2}.
  \]
  Outward rounding gives $\lfloor Qs^2\rfloor/Q-\widehat E>0.3$.
  The triple
  \[
    \left(\frac{a}{1+\widehat E}, \frac{b}{\lfloor Qs^2\rfloor/Q-\widehat E},
      \frac{2\widehat Z^2}{(\lfloor Q(2\pi_-/n)^2\rfloor/Q)(1
      +\sqrt{1+\widehat Z^2})}\right)
  \]
  lies in the corresponding product of intervals above.
  In the last component replace $\widehat Z^2$ by $Q^{-1}\lceil Q\widehat Z^2\rceil$
  before rounding the denominator down, since $z/(1+\sqrt{1+z})=\sqrt{1+z}-1$
  increases for $z\ge0$.
\end{proof}

\begin{proof}[Proof of Numerical Lemma~\ref{sector-spectral-intervals}]
  Let $h=\pi/n$, $z(t)=\cos h+it\sin h$ for $-1\le t\le1$, and
  \[
    M_{jk}=\frac12\int_{-1}^1z(t)^j\overline{z(t)}^k\,dt=M_{kj}
    \quad(j,k\ge0),\qquad D_l=\frac{\sin(lh)}{\sin h}\quad(l\in\mathbb Z).
  \]
  Integrating $(z^k\bar z^{j+1})'$ on $[-1,1]$ gives
  \begin{equation}
    \label{regular-sector-rounding}
    M_{j0}=\frac{D_{j+1}}{j+1},\qquad
    M_{jk}=\frac{2k\cos h\,M_{j,k-1}-D_{k-1-j}}{j+k+1}\quad(1\le k\le j).
  \end{equation}
  Take the smaller root $\lambda$ of the Ritz equation on
  $\mathrm{span}\{z,\bar z^{n-1}\}$,
  \[
    (\cos h-\lambda M_{11})\big((n-1)\cos h\,M_{n-2,n-2}-\lambda M_{n-1,n-1}\big)
      -\lambda^2 M_{n0}^2=0,
  \]
  and set $\mu=10^{-12}\lfloor10^{12}\lambda+1/2\rfloor$.
  On the vertical side $z(t)$ the residuals $(\partial_x-\mu)z^{1+nj}$ for $j\ge0$ and
  $(\partial_x-\mu)\bar z^{-1-nj}$ for $j<0$ are
  \[
    F_j(t)=
    \begin{cases}
      (1+nj-\mu z(t))z(t)^{nj}, & 0\le j\le16, \\
      (-1-nj-\mu\overline{z(t)})\overline{z(t)}^{-2-nj}, & -16\le j<0 .
    \end{cases}
  \]
  With the integrals real by reflection, solve
  \[
    \begin{gathered}
      H_{ij}=\left\lfloor\frac{10^{36}}2\int_{-1}^1 F_i(t)\overline{F_j(t)}\,dt
        +\frac12\right\rfloor
        \quad(-16\le i,j\le16), \\
      10^{20}c_i+\sum_{1\le|j|\le16}H_{ij}c_j=-H_{i0}\quad(1\le|i|\le16).
    \end{gathered}
  \]
  The unrounded matrix is a Gram matrix.
  For any real $v_j$, $1\le|j|\le16$,
  \[
    \sum_{1\le|i|\le16}\sum_{1\le|j|\le16}v_iH_{ij}v_j
    \ge-\tfrac12\left(\sum_{1\le|j|\le16}|v_j|\right)^2
    \ge-16\sum_{1\le|j|\le16}v_j^2 .
  \]
  The system thus has a unique solution.
  Round each $c_j$ to $10^{-12}\lfloor10^{12}c_j+1/2\rfloor$, set $c_0=1$, and use the
  polynomial
  \[
    f(z)=\sum_{j=0}^{16}c_jz^{nj+1}+\sum_{j=1}^{16}c_{-j}\bar z^{nj-1}.
  \]
  For the side $E$ and its sector $S$, radial integration gives
  \[
    \frac1{2\sin h}\int_E z^j\bar z^k\,ds=M_{jk},\qquad
    \frac1{\sin h\cos h}\int_S z^j\bar z^k\,dx\,dy=\frac{2M_{jk}}{j+k+2}
    \quad(j,k\ge0).
  \]
  Substitute
  \[
    (\partial_x f,\partial_y f)=(1,i)\sum_{j=0}^{16}(nj+1)c_jz^{nj}
      +(1,-i)\sum_{j=1}^{16}(nj-1)c_{-j}\bar z^{nj-2}
  \]
  in the sector moments of $f$.
  Use $\sin h=\sqrt{1-\cos^2h}$ and
  \[
    \sum_{i=0}^{31}\frac{(-1)^i(\pi_+/n)^{2i}}{(2i)!}
    \le\cos h\le\sum_{i=0}^{30}\frac{(-1)^i(\pi_-/n)^{2i}}{(2i)!}.
  \]
  In \eqref{regular-sector-rounding} use $D_0=0$, $D_1=1$,
  $D_{l+1}=2\cos h\,D_l-D_{l-1}$ for $1\le l\le n-2$, $D_n=0$, $D_{-l}=-D_l$ and
  $D_{n+l}=-D_l$.
  Outward rounding on $10^{-60}\mathbb Z$ fixes all $H_{ij}$ and $10^{12}\mu$.
  Substitution in \eqref{eq:sector-trial-residual} and its moment bounds gives
  \eqref{eq:finite-scaled-moments}.
  More explicitly, for $5\le n\le20$,
  $0<\mu<g=(n-1)\cos h/(\sqrt{1+\tan^2 h/4}+\tan h/2)^2$ and
  \[
    \begin{aligned}
      \frac{\int_E|\partial_x f-\mu f|^2\,ds}{(g-\mu)^2}
        &<(0.0036)^2\int_E|f|^2\,ds, \\
      \frac{g\int_E|\partial_x f-\mu f|^2\,ds}{(g-\mu)^2}
        &<(0.0044)^2\int_S\big(|\partial_x f|^2+|\partial_y f|^2\big)\,dx\,dy .
    \end{aligned}
  \]

  For $21\le n\le314$ one has $h<0.15$, $\tan h<0.152$, $\cos h>0.9887$ and
  $\sin h/h>0.9962$.
  Thus $\gamma_n>0.8n$.
  Since $1-3.044/n^3<s_n\le q_n\le1$,
  \[
    0<\frac{1-S_n/q_n}{1-S_n/\gamma_n}<\frac{3.24}{n^3},\qquad
    \frac{q_n-S_n}{\gamma_n-S_n}<\frac{4.05}{n^4}<
    \frac{\sin^2 h}{3\cos^2 h+\sin^2 h}<\frac12 .
  \]
\end{proof}

\begin{proof}[Proof of Numerical Lemma~\ref{certbern}]
  Substitution of \eqref{sect-matrix} into \eqref{eq:sector-polynomial-defect} gives
  \begin{align*}
    P(r)={} & 27(T/2+mV)r^3
      \left(r^2+\frac{d(1+pr^2)}{4p}\right) \\
    & \quad\times\left[b\left(1+\frac{d(1+pr^2)}4\right)
      +(1-b)\left(r^2+\frac{d(1+pr^2)}{4p}\right)\right]^2 \\
    & -4T\bigg[(3/2-B-C)r^2+B\left(1+\frac{d(1+pr^2)}4\right)^2 \\
    & \qquad+C\left(r^2+\frac{d(1+pr^2)}{4p}\right)^2
      +\frac{Kd(d+4)(1+pr^2)^2}{16p}\bigg]^3 .
  \end{align*}
  In \eqref{eq:explicit-bernstein-coefficients} use Taylor's formula
  \[
    P(1+\delta w)=\sum_{i=0}^{\deg P}\frac{\delta^iP^{(i)}(1)}{i!}\,w^i .
  \]
  For $p=\tan^2(\pi/n)$,
  \[
    \begin{aligned}
      \sum_{i=0}^{15}\frac{(-1)^i(\pi_-/n)^{2i+1}}{(2i+1)!}
        & \le\sin(\pi/n)\le
        \sum_{i=0}^{14}\frac{(-1)^i(\pi_+/n)^{2i+1}}{(2i+1)!}, \\
      \sum_{i=0}^{15}\frac{(-1)^i(\pi_+/n)^{2i}}{(2i)!}
        & \le\cos(\pi/n)\le
        \sum_{i=0}^{14}\frac{(-1)^i(\pi_-/n)^{2i}}{(2i)!}.
    \end{aligned}
  \]
  Substitution at the $16$ vertices of \eqref{eq:finite-scaled-moments} for
  $5\le n\le20$ and the $12$ vertices of \eqref{eq:coefficient-vertices} for
  $21\le n\le314$ gives the bounds, with outward rounding on $10^{-48}\mathbb Z$.
\end{proof}

\begin{numericallemma}[Five pentagonal perturbations]
  \label{certpentagon}
  On the regular pentagon $R_5$ with vertices $e^{2\pi ij/5}$, $0\le j<5$, take the
  harmonic polynomial
  \[
    \begin{aligned}
      10^5 f(z)={} & -78\bar z^{39}+100\bar z^{34}-64\bar z^{29}-3\bar z^{24}
        -30\bar z^{19}-70\bar z^{14}-254\bar z^9-6399\bar z^4 \\
      & -45726z+1390z^6+158z^{11}+50z^{16}+16z^{21}
        +50z^{26}-78z^{31}+72z^{36}.
    \end{aligned}
  \]
  Write $z=x+iy$.
  On the open side from $e^{2\pi ij/5}$ to $e^{2\pi i(j+1)/5}$ the outward unit normal
  is
  \[
    (\nu_1,\nu_2)=\left(\cos\frac{(2j+1)\pi}{5},\sin\frac{(2j+1)\pi}{5}\right)
    \qquad(0\le j<5).
  \]
  With $ds$ denoting arclength,
  \begin{equation}
    \label{pen-resid}
    \begin{aligned}
      \left|\int_{\partial R_5}|f|^2\,ds-1\right| &<10^{-3}, \\
      \left|\int_{R_5}\big(|\partial_x f|^2+|\partial_y f|^2\big)\,dx\,dy
        -1.017\right| &<10^{-4}, \\
      \int_{\partial R_5}|\nu_1\partial_x f+\nu_2\partial_y f
        -1.017f|^2\,ds &<(0.003)^2, \\
      \int_{R_5}\big(|\partial_x f+0.45726|^2+|\partial_y f+0.45726i|^2\big)\,dx\,dy
        &<0.0416 .
    \end{aligned}
  \end{equation}
  For the displacement $W$ and the coefficients in \eqref{pen-hess}, put
  \begin{align*}
    \tfrac12 Q_f(W)={} & \int_{R_5}\big(\langle H\nabla f,\nabla f\rangle+
      2\operatorname{Re}\big(\langle G\nabla f,\nabla Z\rangle\big)
        +|\nabla Z|^2\big)\,dx
      +\tfrac{1.017}{2}\kappa \\
    & -1.017\int_{\partial R_5}\big(\beta|f|^2+2\alpha\operatorname{Re}(f\overline Z)
      +|Z|^2\big)\,ds
      +\frac{1.017}{L(R_5)}\left|\int_{\partial R_5}(\alpha f+Z)\,ds\right|^2 .
  \end{align*}
  For the modes $w_2,w_3,w_4$ defined in the proof of Theorem~\ref{pentlocaltheorem},
  \begin{equation}
    \label{pen-five}
    \begin{gathered}
      Q_f(w_2)<-0.25,\quad Q_f(w_3)<-0.78,\quad Q_f(w_4)<-0.79, \\
      \left|Q_f\!\left(\tfrac{w_3+w_4}{\sqrt2}\right)
        -Q_f\!\left(\tfrac{w_3-w_4}{\sqrt2}\right)\right|<0.75 .
    \end{gathered}
  \end{equation}
  For each of these displacements and each $0\le j<5$,
  \begin{equation}
    \label{pen-coeff}
    \begin{gathered}
      \max_{|v|=1}|H_jv|<1.05,\quad \max_{|v|=1}|G_jv|<1.27,\quad
        |\alpha_j|<0.83,\quad |\beta_j|<0.32, \\
      |\kappa|<0.21,\qquad
        \int_{R_5}|\nabla Z|^2\,dx\,dy<(0.53)^2,\qquad
        \int_{\partial R_5}|Z|^2\,ds<(0.44)^2 .
    \end{gathered}
  \end{equation}
\end{numericallemma}
\begin{proof}
  Put $v_j=\omega^j$, $0\le j\le5$, where $\omega=(\sqrt5-1+i\sqrt{10+2\sqrt5})/4$.
  On side $j$ the length is $2\sin(\pi/5)$ and the outward unit normal is
  $\nu_j=e^{(2j+1)\pi i/5}$.
  For a polynomial $h$ of total degree $D$ in $z,\bar z$, expand
  \[
    h\big(r((1-t)v_j+tv_{j+1})\big)
    =\sum_{m=0}^D\sum_{q=0}^m h_{mq}r^mt^q .
  \]
  Integrating on $0\le r,t\le1$ gives
  \[
    \begin{aligned}
      \int_{S_j}h\,dx
        &=\sin(2\pi/5)\sum_{m=0}^D\sum_{q=0}^m\frac{h_{mq}}{(m+2)(q+1)}, \\
      \int_{[v_j,v_{j+1}]}h\,ds
        &=2\sin(\pi/5)\sum_{m=0}^D\sum_{q=0}^m\frac{h_{mq}}{q+1}.
    \end{aligned}
  \]
  For \eqref{pen-resid} substitute $|f|^2$, $|f_x|^2+|f_y|^2$,
  $|\operatorname{Re}(\nu_j)\,\partial_x f+\operatorname{Im}(\nu_j)\,\partial_y f
  -1.017f|^2$ and $|\partial_x f+0.45726|^2+|\partial_y f+0.45726i|^2$ in these
  formulas and sum over $0\le j\le4$.
  For \eqref{pen-five} the interpolation formula \eqref{pen-interpolant} gives
  \begin{align*}
    Z\big(r((1-t)v_j+tv_{j+1})\big) &=-0.45726\,r\big((1-t)W_j+tW_{j+1}\big), \\
    \partial_{x_k} Z &=-0.45726\big((B_j)_{1k}+i(B_j)_{2k}\big)\quad(k=1,2).
  \end{align*}
  Substitution in $Q_f$ gives \eqref{pen-five}.
  For \eqref{pen-coeff}, use
  \[
    2\max_{|v|=1}|Mv|=|M_{11}+M_{22}|+\sqrt{(M_{11}-M_{22})^2+4M_{12}^2}
    \qquad(M=G_j,H_j).
  \]
  All monomial coefficients and square roots are enclosed on $10^{-80}\mathbb Z$.
\end{proof}

\begin{numericallemma}[Third variation at the pentagon]
  \label{pent-taylor}
  Under the normalization of Theorem~\ref{pentlocaltheorem}, assume
  $\sum_{j=0}^4 W_j=\sum_{j=0}^4 W_j\bar v_j=0$ and $\sum_{j=0}^4|W_j|^2=1$.
  For $0\le t\le1/100$,
  \[
    \frac{|D'''(t)|}{6}<6.3,\qquad
    0\le\int_{\partial P_t}|u_t|^2\,ds-
      \frac{1}{L(P_t)}\left|\int_{\partial P_t}u_t\,ds\right|^2<1.025 .
  \]
\end{numericallemma}
\begin{proof}
  By \eqref{pen-interpolant} and Cauchy--Schwarz,
  \[
    |\xi_j|\le\max_{|v|=1}|B_jv|
    =\sup_{(s,t)\in\mathbb R^2\setminus\{(0,0)\}}\frac{|sW_j+tW_{j+1}|}
      {|sv_j+tv_{j+1}|}
    \le\sqrt{\frac{|W_j|^2+|W_{j+1}|^2}{1-\cos(2\pi/5)}}<1.21
    \qquad(0\le j<5).
  \]
  Thus $x\mapsto x+tV(x)$ is bi-Lipschitz for $0\le t\le0.01$.
  Fourier orthogonality gives $\int_{\partial R_5}Z\,ds=0$ and
  \begin{align*}
    \kappa &=\sum_{k=2}^4\frac{\sin(2\pi k/5)}{5\sin(2\pi/5)}|z_k|^2,
      & |\kappa| & \le 1/5, \\
    \int_{\partial R_5}|Z|^2\,ds &=(0.45726)^2 L(R_5)\sum_{k=2}^4
      \frac{2+\cos(2\pi k/5)}{15}|z_k|^2
      <(0.44)^2 .
  \end{align*}
  In the derivative estimates below, $g_i$ denotes an upper bound for
  $\sup_{0\le t\le0.01}|g^{(i)}(t)|/i!$, $0\le i\le3$.
  Put $e(x)=U(x)+0.45726x$, identifying $x\in\mathbb R^2$ with its complex coordinate.
  By \eqref{pen-resid} and the projection in Theorem~\ref{pentlocaltheorem},
  $\int_{R_5}(|\partial_x e|^2+|\partial_y e|^2)\,dx\,dy<0.045$.
  On $S_j$, $0\le j\le4$, set $J=I_2+tB_j$ and
  $F=(\operatorname{adj}J)(\operatorname{adj}J)^{\mathsf T}/\det J$, where
  $\operatorname{adj}J=(\det J)J^{-1}$.
  Since $u_t(x+tV(x))=-0.45726(x+tV(x))+e(x)$,
  \begin{align*}
    E(t)=\int_{P_t}|\nabla u_t|^2\,dx
      &=2(0.45726)^2|R_5|(1+\kappa t^2) \\
    & \quad-2(0.45726)\int_{R_5}\operatorname{tr}(De\,\operatorname{adj}J)\,dx
      +\int_{R_5}\langle F\nabla e,\nabla e\rangle\,dx .
  \end{align*}
  The adjugate is affine in $t$ with $\max_{|v|=1}|(\operatorname{adj}J)'v|<1.21$.
  The eigenvalues $\lambda_1,\lambda_2$ of $B_j$ satisfy
  $|\lambda_1|,|\lambda_2|<1.21$ and $\det J=(1+t\lambda_1)(1+t\lambda_2)$.
  Hence
  \[
    \frac{\max_{|v|=1}|F^{(i)}(t)v|}{i!}
    \le\sum_{j=0}^{\min(i,2)}\binom2j
    \frac{(1.0121)^{2-j}(i-j+1)(1.21)^i}{(0.9879)^{i-j+2}},
    \qquad 0\le i\le3 .
  \]
  Since $|R_5|=\tfrac52\sin(2\pi/5)<2.378$, Cauchy--Schwarz gives, componentwise,
  \[
    \begin{aligned}
      \left(\sup_{0\le t\le0.01}\frac{|E^{(i)}(t)|}{i!}\right)_{i=0}^3 &
        \le 2(0.45726)^2|R_5|(1.00002,0.004,1/5,0) \\
      & \quad+2(0.45726)\sqrt{2|R_5|\cdot0.045}\,(1.0121,1.21,0,0) \\
      & \quad+0.045(1.050,5.081,12.372,22.685) \\
      &<(1.5,0.75,0.8,1.03).
    \end{aligned}
  \]
  As $|\kappa|\le1/5$, Leibniz' rule gives
  \[
    \begin{aligned}
      \frac{1}{6}\left|\frac{d^3}{dt^3}\left(\sqrt{1+\kappa t^2}\,E(t)\right)\right|
        & \le1.03\sqrt{1.00002}+\frac{0.008}{5\sqrt{0.99998}} \\
      & \quad+\frac{0.75}{10(0.99998)^{3/2}}+
        \frac{0.015}{50(0.99998)^{5/2}}<1.11 .
    \end{aligned}
  \]

  On side $j$, $0\le j\le4$, put $w(t)=|1+t\xi_j|$.
  The identities $ww''=|\xi_j|^2-(w')^2$, $ww'''=-3w'w''$ and the bound $|w'|\le1.21$
  give
  \[
    (w_0,w_1,w_2,w_3)=
    \left(1.0121,1.21,\frac{1.21^2}{2(0.9879)},\frac{1.21^3}{2(0.9879)^2}\right).
  \]
  Since $\int_{\partial R_5}|U|^2\,ds<(1.001)^2$ and
  $\int_{\partial R_5}|Z|^2\,ds<(0.44)^2$,
  \[
    \begin{aligned}
      \int_{\partial R_5}|U+tZ|^2\,ds & \le(1.001+0.44t)^2<(1.006)^2, \\
      \frac1{6}\left|\frac{d^3}{dt^3}\int_{\partial P_t}|u_t|^2\,ds\right|
        & \le w_3(1.006)^2+2w_2(1.006)(0.44)+w_1(0.44)^2<1.81 .
    \end{aligned}
  \]
  The centering term equals $|M(t)|^2/\ell(t)$, where
  \[
    M(t)=L(R_5)^{-1/2}\int_{\partial R_5}(U+tZ)w\,ds,\qquad
    \ell(t)=L(R_5)^{-1}L(P_t)=L(R_5)^{-1}\int_{\partial R_5}w\,ds .
  \]
  Since $M(0)=0$ by rotation and
  $\ell'(0)=\tfrac15\operatorname{Re}\left(\sum_{j=0}^4 W_j\bar v_j\right)=0$,
  convexity gives $\ell(t)\ge1$.
  By $\int_{\partial R_5}Z\,ds=0$, $|M'(0)|\le1.001\cdot1.21$ and
  \[
    M_2=1.006w_2+0.44w_1,\quad M_3=1.006w_3+0.44w_2,\quad
    M_1=1.001\cdot1.21+0.02M_2,\quad M_0=0.01M_1 .
  \]
  With $(M_0,M_1,M_2,M_3)\le(0.01237,1.23677,1.27787,1.23912)$ and
  $|\ell'(t)|\le0.02w_2$, Leibniz' rule gives
  \[
    \begin{aligned}
      \frac1{6}\left|\frac{d^3}{dt^3}\frac{|M|^2}{\ell}\right|
        & \le 2(M_0M_3+M_1M_2)+0.02w_2(2M_0M_2+M_1^2) \\
      & \quad +2M_0M_1((0.02w_2)^2+w_2)
        +M_0^2((0.02w_2)^3+0.04w_2^2+w_3)<3.26 .
    \end{aligned}
  \]
  Thus $|D'''(t)|/6<1.11+1.021(1.81+3.26)<6.3$.
  By Cauchy--Schwarz the centered boundary mass lies between zero and
  $w_0(1.006)^2<1.025$.
\end{proof}

\bigskip
\noindent
{\bf Acknowledgements.} The authors used AI to assist the long but finite computations.


\begin{thebibliography}{99}

\bibitem{AKO}
  E. Akhmetgaliyev, C.-Y. Kao and B. Osting,
  Computational methods for extremal Steklov problems,
  {\em SIAM J. Control Optim.} \textbf{55} (2017), 1226--1240.

\bibitem{AIM}
  American Institute of Mathematics,
  \emph{Steklov eigenproblems},
  workshop problem list, Section 4, Problem 4.6 (\'E. Oudet), 2018.

\bibitem{BBF}
  B. Bogosel, D. Bucur and I. Fragal\`a,
  Proof of the polygonal Saint-Venant inequality conjectured by
    P\'olya--Szeg\H{o},
  arXiv:2611.24573, 2026.

\bibitem{BDR}
  L. Brasco, G. De Philippis and B. Ruffini,
  Spectral optimization for the Stekloff--Laplacian: the stability issue,
  {\em J. Funct. Anal.} \textbf{262} (2012), 4675--4710.

\bibitem{Brock}
  F. Brock,
  An isoperimetric inequality for eigenvalues of the Stekloff problem,
  {\em Z. Angew. Math. Mech.} \textbf{81} (2001), 69--71.

\bibitem{BFNT}
  D. Bucur, V. Ferone, C. Nitsch and C. Trombetti,
  Weinstock inequality in higher dimensions,
  {\em J. Differential Geom.} \textbf{118} (2021), 1--21.

\bibitem{CGHL}
  Z. Cheng, C. Gui, Y. Hu and Q. Li,
  The P\'olya--Szeg\H{o} conjecture for convex polygons with many sides,
  arXiv:2609.18500, 2026.

\bibitem{CGHLY}
  Z. Cheng, C. Gui, Y. Hu, Q. Li and R. Yao,
  Global strict monotonicity and asymptotics of the first Steklov eigenvalue for
    regular polygons,
  arXiv:2603.25116v2, 2026.

\bibitem{DGP}
  J. Dahne, J. G\'omez-Serrano and J. Pech-Alberich,
  Monotonicity of the first Dirichlet eigenvalue of regular polygons,
  arXiv:2601.16285, 2026.

\bibitem{D}
  J. Davies,
  The equilateral triangle maximises the first Steklov eigenvalue among
    triangles,
  preprint, Hexagon 2611.00629v1 (\nolinkurl{https://hexagonmath.org/2611.00629v1}),
    2026.

\bibitem{DNS}
  S. Dom\'inguez, N. Nigam and B. Shahriari,
  A combined finite element and Bayesian optimization framework for shape
    optimization in spectral geometry,
  {\em Comput. Math. Appl.} \textbf{74} (2017), 2874--2896.

\bibitem{DYP}
  L. Dong, S. Yang and J. Pan,
  The square maximizes the first nonzero Steklov eigenvalue among convex
    quadrilaterals,
  preprint, Hexagon 2611.00777v1 (\nolinkurl{https://hexagonmath.org/2611.00777v1}),
    2026.

\bibitem{DrTr}
  T. A. Driscoll and L. N. Trefethen,
  \emph{Schwarz--Christoffel Mapping},
  Cambridge University Press, Cambridge, 2002.

\bibitem{DGMRV}
  E. B. Dryden, C. S. Gordon, J. Moreno, J. Rowlett and C. Villegas-Blas,
  The Steklov spectrum of convex polygonal domains I: spectral finiteness,
  {\em J. Geom. Anal.} \textbf{35} (2025), no.~3, Paper No.~91.

\bibitem{Farouki}
  R. T. Farouki,
  The Bernstein polynomial basis: a centennial retrospective,
  {\em Comput. Aided Geom. Design} \textbf{29} (2012), 379--419.

\bibitem{FL}
  P. Freitas and R. S. Laugesen,
  From Steklov to Neumann and beyond, via Robin: the Szeg\H{o} way,
  {\em Canad. J. Math.} \textbf{72} (2020), 1024--1043.

\bibitem{GLPS}
  A. Girouard, J. Lagac\'e, I. Polterovich and A. Savo,
  The Steklov spectrum of cuboids,
  {\em Mathematika} \textbf{65} (2019), 272--310.

\bibitem{GP}
  A. Girouard and I. Polterovich,
  Spectral geometry of the Steklov problem (survey article),
  {\em J. Spectr. Theory} \textbf{7} (2017), 321--359.

\bibitem{GHL}
  C. Gui, Y. Hu and Q. Li,
  A quantitative P\'olya--Szeg\H{o} theorem for tangential polygons,
  arXiv:2607.28768, 2026.

\bibitem{GHLY}
  C. Gui, Y. Hu, Q. Li and H. Yang,
  A quantitative Robin P\'olya--Szeg\H{o} theorem for tangential polygons,
  arXiv:2610.21417, 2026.

\bibitem{GHLZ}
  C. Gui, Y. Hu, Q. Li and C. Zhang,
  Mixed torsion on right triangles and the P\'olya--Szeg\H{o} monotonicity
    problem for regular polygons,
  arXiv:2606.13448, 2026.

\bibitem{Henrot}
  A. Henrot (ed.),
    {Shape Optimization and Spectral Theory},
  De Gruyter Open, Warsaw, 2017.

\bibitem{HPS}
  J. Hersch, L. E. Payne and M. M. Schiffer,
  Some inequalities for Stekloff eigenvalues,
  {\em Arch. Rational Mech. Anal.} \textbf{57} (1974), 99--114.

\bibitem{IN}
  E. Indrei and L. Nurbekyan,
  On the stability of the polygonal isoperimetric inequality,
  {\em Adv. Math.} \textbf{276} (2015), 62--86.

\bibitem{LS}
  R. S. Laugesen and B. A. Siudeja,
  Maximizing Neumann fundamental tones of triangles,
  {\em J. Math. Phys.} \textbf{50} (2009), 112903.

\bibitem{LPPS}
  M. Levitin, L. Parnovski, I. Polterovich and D. A. Sher,
  Sloshing, Steklov and corners: asymptotics of Steklov eigenvalues for
    curvilinear polygons,
  {\em Proc. Lond. Math. Soc.} (3) \textbf{125} (2022), 359--487.

\bibitem{YLi}
  Y. Li,
  A proof of the polygonal Faber--Krahn inequality,
  arXiv:2611.27968, 2026.

\bibitem{LX}
  S. Luo and B. Xie,
  The regular $N$-gon maximizes first positive Neumann eigenvalue among simple
    polygons with at most $N$ sides and fixed area,
  arXiv:2611.26729, 2026.

\bibitem{MKC}
  R. E. Moore, R. B. Kearfott and M. J. Cloud,
  \emph{Introduction to Interval Analysis},
  SIAM, Philadelphia, 2009.

\bibitem{Nigam}
  N. Nigam,
  Polygonal Weinstock-type conjecture,
  Problem~8 in \emph{Open problems from the miniconference on sharp eigenvalue
    estimates for partial differential operators}, April 2020, p.~2,
  problem list (\nolinkurl{https://publish.illinois.edu/eigenvalues2020/files/2020/04/Open-problems.pdf}).

\bibitem{PS}
  G. P\'olya and G. Szeg\H{o},
  \emph{Isoperimetric Inequalities in Mathematical Physics},
  Annals of Mathematics Studies, vol.~27, Princeton University Press, Princeton,
    NJ, 1951.

\bibitem{Wang}
  H. Wang,
  The P\'olya--Szeg\H{o} inequalities for polygons,
  arXiv:2611.22305, 2026.

\bibitem{W}
  R. Weinstock,
  Inequalities for a classical eigenvalue problem,
  {\em J. Rational Mech. Anal.} \textbf{3} (1954), 745--753.
\end{thebibliography}
\end{document}